\documentclass[11pt]{article}
\usepackage{graphicx, color}
\usepackage{graphicx}

\def\la{\big\langle}
\def\ra{\big\rangle}

\def\ds{\displaystyle}
\def\forall{\hbox{for all}~}
\def\L{{\bf L}}

\def\bfv{{\bf v}}

\def\bfw{{\bf w}}

\def\bfn{{\bf n}}

\def\ve{\varepsilon}

\def\R{I\!\!R}

\def\div{\hbox{div}}

\def\implies{\Longrightarrow}
\def\vp{\varphi}

\def\vsk{\vskip 4em}
\def\v{\vskip 1em}
\def\O{{\cal O}}

\def\begi{\begin{itemize}}
\def\endi{\end{itemize}}
\def\C{{\cal C}}

\def\ov{\overline}
\def\Tilde{\widetilde}

\def\bega{\begin{array}}
\def\enda{\end{array}}
\def\meas{\hbox{meas}}

\def\bel{\begin{equation}\label}
\def\eeq{\end{equation}}
\def\sqr#1#2{\vbox{\hrule height .#2pt
\hbox{\vrule width .#2pt height #1pt \kern #1pt
\vrule width .#2pt}\hrule height .#2pt }}
\def\square{\sqr74}
\def\endproof{\hphantom{MM}\hfill\llap{$\square$}\goodbreak}

\begin{document}
\title{\bf  Approximation of Sweeping Processes and Controllability 
for a Set Valued Evolution}

\author{Alberto Bressan$^{(1)}$, Marco Mazzola$^{(2)}$, and Khai T. Nguyen$^{(3)}$\\ 
\\
 {\small $^{(1)}$ Department of Mathematics, Penn State University, }\\  
 {\small $^{(2)}$ IMJ-PRG, CNRS, Sorbonne Universit\'e,}\\
  {\small $^{(3)}$ Department of Mathematics, North Carolina State University. }\\ 
 \\  {\small e-mails: ~axb62@psu.edu, ~marco.mazzola@imj-prg.fr,  ~ khai@math.ncsu.edu}
 }
\maketitle
\begin{abstract}
We consider  a controlled evolution problem for a set $\Omega(t)\in\R^d$,
originally motivated by a model where a dog controls a flock of sheep.
Necessary conditions and sufficient conditions are given, in order that the 
evolution be completely controllable.  Similar techniques are then applied
to the approximation of a sweeping process.
Under suitable assumptions, we prove that there exists a control
function such that the corresponding evolution of the set $\Omega(t)$ is arbitrarily close to the one 
determined by the sweeping process.
\end{abstract}

\vsk

\section{Introduction}
\label{sec:1}
\setcounter{equation}{0}
In this paper we consider a controllability problem for the evolution of a set 
$\Omega(t)\subset \R^d$.  This was originally motivated by
the model introduced in \cite{BZ}, describing the evolution of  a flock
of sheep, who tend to scatter around but also react to the presence of a dog.
The region $\Omega(t)\subset \R^2$ occupied by the sheep is 
described as the reachable set for a
differential inclusion, while the position of the dog is regarded as a control function. 
As in \cite{BZ}, we consider a ``scare function"  $\vp=\vp(r)>0$, describing 
the speed at which sheep run away from the dog, depending on the distance $r$. 
Further results and extensions can be found in \cite{CP1, CP2}. For more general models
of crowd dynamics we refer to \cite{BD}.  A general theory of evolution
problems in metric spaces, also describing the evolution of a set, was developed in \cite{Aubin, L}.

In the following we consider the evolution of a set in $\R^d$, and assume
\begi
\item[{\bf (A1)}] {\it The function $r\mapsto \vp(r) $ is continuously differentiable 
for $r>0$,  and satisfies}
\bel{vplim}\vp'~<~0,
\qquad\quad
\lim_{r\to 0+}\vp(r)~=~+\infty,\qquad\quad \lim_{r\to +\infty}\vp(r)~=~ 0.\eeq
\endi

Given a function $t\mapsto \xi(t)\in \R^d$ describing the position of a repelling agent,
we define the velocity field
\bel{G}
\bfv(x,\xi)~\doteq~
\vp(|x-\xi|)\,{x-\xi\over |x-\xi|}\,.\eeq
For a given initial set $\Omega_0$, we denote by $\Omega^\xi(t)$ the 
set reached by trajectories of 
\bel{ODE1} \dot x~\in~\bfv(x, \xi(t))\,,\qquad\qquad x(0) ~\in \Omega_0\,.\eeq
In other words, for any $t\geq 0$,
\bel{Stdef}\bega{rl}
\Omega^\xi(t)~\doteq~\Big\{x(t)\,;  &~x(0)\in \Omega_0\,,\qquad x(\cdot)~~\hbox{is absolutely continuous},\cr
&\qquad\qquad ~~~\dot x(\tau)~=~
\bfv\left(x(\tau),\,\xi(\tau)\right) ~~\hbox{for a.e.}~\tau\in [0,t]\Big\}.\enda
\eeq
Throughout  the following we write $\partial \Omega$, $\ov\Omega$, and $int\,\Omega$, 
 for the boundary, the closure, and the interior of a set $\Omega\subset\R^d$, 
respectively. 
By $B(\Omega,r)$ we denote the open neighborhood of radius $r$ around the set $\Omega$,
while 
$d_H$ denotes Hausdorff distance \cite{AC}.

To avoid any difficulty about uniqueness of solutions of (\ref{G})-(\ref{ODE1}), we assume
that the control $\xi(\cdot)$ is chosen so that 
\bel{Ad}
\inf_{t\in [0,\tau]}~d\bigl(\xi(t),\Omega^\xi(t)\bigr)>0\qquad\qquad\forall 0\leq\tau<T.
\eeq
We wish to understand how the function $\vp$ affects the controllability properties of the evolution equation
(\ref{ODE1}).  Roughly speaking, given an initial set $\Omega_0$ and a terminal set $\Omega_1$, we seek  a control $\xi(\cdot)$  such that, at the terminal time $T$, the 
set $\Omega^\xi(T)$ in (\ref{Stdef}) is arbitrary close to $\Omega_1$.
\v
{\bf Definition 1.} {\it
We say that the set-valued evolution (\ref{ODE1}) satisfies the  {\bf Global Approximate
Confinement} property  {\bf (GAC)} if the following holds.
Let  $\Omega_0,\Omega_1 \subset\R^d$ be any two compact domains,
with $\Omega_1$ simply connected and such that $\Omega_1\subset int\,\Omega_0$.
Then for any $T,\ve>0$, there exists a Lipschitz continuous 
 control $\xi:[0,T]\mapsto\R^d$ satisfying (\ref{Ad}) and such that
the set (\ref{Stdef})
satisfies
\bel{DHT}  \Omega_1~\subseteq~\Omega^\xi(T)~\subseteq~B(\Omega_1,\ve)\eeq
If there exists a locally Lipschitz control $\xi:[0,T[\mapsto \R^d$ 
satisfying (\ref{Ad}) and such that 
\bel{D0}
\Omega^\xi(T)~=~ \Omega_1\,,\eeq
we then say that the set-valued evolution (\ref{ODE1})  satisfies the {\bf Global Exact 
Confinement} property {\bf (GEC)}. }
\v
The primary goal of this paper
is to find conditions which are necessary, or sufficient, to achieve the {\bf (GAC)} or 
{\bf (GEC)} properties.  Indeed, we will show that these properties are determined
by the asymptotic behavior of the function $\vp$ as $r\to 0+$.  
Our first main result is
\v
{\bf Theorem 1 (necessary condition).}  {\it Let $\vp$ satisfy {\bf (A1)}. If the {\bf (GAC)}
property holds, then  the function $\vp$ must satisfy}
\bel{EC3}
\int_{0}^{1}r^{d-2}\varphi(r)\,dr~=~+\infty\,.
\eeq

To state a sufficient condition, we introduce the assumption
\begi
\item[{\bf (A2)}] {\it For every $\kappa>0$ one has}
\bel{lim0}
\lim_{r\to 0+}~{r^{d/2}\cdot \vp(\kappa r^{1/2})+1\over r^d\cdot \vp(r)}~=~0\,.
\eeq
\endi
\v
{\bf Theorem 2 (sufficient condition).}  {\it If the function $\vp$ satisfies 
{\bf (A1)-(A2)}, then the {\bf (GAC)} and {\bf (GEC)} properties hold.}
\v
This result applies, in particular, to the function $\vp(r) = r^{-\beta}$ for any $\beta> d$.

\v
A proof of Theorem 1 will be given in Section~2, while Theorem~2 will be proved in 
Section~3.   

The controllability of the set-valued evolution  (\ref{Stdef})
is closely related to a result on the approximation of a sweeping process. 
Indeed,
let $t\mapsto V(t)$ describe a moving set in $\R^d$.  We assume that 
each $V(t)$ is a compact set with nonempty interior and
smooth boundary, smoothly depending on time. More precisely:
\bel{Vdef}
V(t)~=~\bigl\{ x\in\R^d\,;~~\psi(t,x)\leq 0\bigr\},\eeq
where $\psi:\R\times\R^d\mapsto\R$ has $\C^2$ regularity and satisfies the nondegeneracy condition
\bel{ppp}\psi(t,x)~=~0\qquad\implies\qquad \nabla_x \psi(t,x)~\not= ~0.\eeq
As usual, we denote by 
$N_{V(t)}(x)$ the outer normal cone to $V(t)$ at a boundary point $x\in \partial V(t)$.    
In the case of an interior point  $x\in int V(t)$, we simply define $N_{V(t)}(x)=\{0\}$.
By the well known theory of sweeping processes \cite{CCPT, CG, CMM, CM, M}, 
for any initial point $x_0\in V(0)$,  the differential
inclusion
\bel{sweep} \dot x(t)~\in~- N_{V(t)} (x(t)),\qquad\quad x(0)= x_0 \eeq
has a unique solution $t\mapsto x(t, x_0)\in V(t)$.
In turn, for a given initial set $\Omega_0\subset V(0)$, one can consider the sets
\bel{OT}
\Omega(t)~\doteq~\bigl\{ x(t, x_0)\,;~~x_0\in \Omega_0\bigr\}.\eeq
A natural question is whether there exists a control $\xi(\cdot)$ such that the corresponding
sets $\Omega^\xi (t)$ in (\ref{Stdef}) remain uniformly close to the sets $\Omega(t)$, 
for all $t\in [0,T]$. It turns out that this is  true, under an assumption
which  slightly strengthens {\bf (A2)}, namely:
\v
\begi
\item[{\bf (A2$'$)}] {\it For some $\beta>1/2$  one has}
\bel{lim01}\lim_{r\to 0+}~{r^{\beta d}\cdot \vp( r^{\beta})+1\over r^d\cdot \vp(r)}~=~0.\eeq
\endi
In the following, $t\mapsto x^\xi(t, x_0)$ denotes the solution to 
\bel{xx0} 
\dot x(t)~=~\bfv(x(t), \xi(t)), \qquad\qquad x(0)=x_0\,,\eeq
with $\bfv$ as in (\ref{G}), while $t\mapsto x(t, x_0)$ is the trajectory of  the 
sweeping process (\ref{sweep}), with the same initial condition.
\v
{\bf Theorem 3 (approximation of a sweeping process).}  {\it Assume that the function $\vp$ satisfies 
{\bf (A1)} and {\bf (A2$'$)}.   As in (\ref{Vdef})-(\ref{ppp}), let $t\mapsto V(t)$ 
be a family of sets with $\C^2$ boundaries.
   Then, for any $T,\ve>0$ there exists
a measurable control $t\mapsto \xi(t)$ such that 
\bel{approx}
\bigl| x^\xi(t, x_0) - x(t, x_0)\bigr|~\leq~\ve\qquad\forall x_0\in V(0), ~~t\in [0,T].\eeq 
}

An immediate consequence of (\ref{approx}) is that, 
for any initial subset $\Omega_0\subseteq V(0)$, the corresponding sets 
$\Omega^\xi(t)$ in (\ref{Stdef}) and $\Omega(t)$ in (\ref{OT}) satisfy
\bel{DHO} d_H\bigl(\Omega^\xi(t),\, \Omega(t)\bigr)~\leq~\ve\qquad\forall 
t\in [0,T].\eeq 

A proof of Theorem~3 will be worked out in Section~4.

\v

\section{Proof of the necessary condition}
\label{sec:2}
\setcounter{equation}{0}

In this section we give a proof of Theorem~1.   The main idea is that, if 
(\ref{EC3}) fails, then for any choice of the control $\xi(\cdot)$ 
the volume of the set $\Omega^\xi$ cannot shrink too fast.  This puts a constraint 
on the sets that can be approximately reached at time $T$.
\v
Let $t\mapsto \xi(t)$ be any admissible control.
Fix any time $t\geq 0$ and let $\bfv= \bfv(\cdot, \xi(t))$ be the vector field in (\ref{G}).
Then, calling $\delta\doteq d(\xi(t), \Omega^\xi(t))$, we compute
\bel{dOX}
\bega{rl}\ds{d\over dt} \meas \Big(\Omega^\xi(t)\Big)&\ds
=~\int_{\Omega^\xi(t)\cap B(\xi(t), 1)} \div 
\,\bfv\, dx+\int_{\Omega^\xi(t)\setminus B(\xi(t), 1)} \div \,\bfv\, dx
\\[4mm]
&\ds \geq~ \omega_{d-1} \int_\delta^1 r^{d-1}\vp'(r) \, dr + \vp'(1)\cdot\meas
\bigl(\Omega^\xi(t)\bigr)\,.\enda
\eeq
Here and in the sequel $\omega_{d-1}$ denotes the $(d-1)$-dimensional measure of the 
surface of the unit ball in $\R^d$. 
An integration by parts yields
\bel{ibp}\bega{rl}\ds
\int_\delta^1 r^{d-1}\vp'(r) \, dr&\ds=~-\int_\delta^1 (d-1)r^{d-2} [\vp(r)-\vp(1)]\, dr
-\delta^{d-1} [\vp(\delta)-\vp(1)]\\[4mm]
&\ds\geq~- \int_\delta^1 (d-1)r^{d-2} \vp(r)\, dr - \delta^{d-1}\vp(\delta).
\enda\eeq
If (\ref{EC3}) fails, then 
\[
\int_{\delta}^1r^{d-2}\vp(r)~dr~\leq~ M~\doteq~\int_{0}^{1}r^{d-2}\vp(r)dr~<~\infty\,.
\]
Since $\vp$ is decreasing, one has
\[
\delta^{d-1}\vp(\delta)~=~{1\over d-1}\cdot \int_{0}^{\delta}s^{d-2}\cdot \vp(\delta)~ds~\leq~{1\over d-1}\cdot \int_{0}^{\delta}s^{d-2}\cdot \vp(s)~ds~\leq~{M\over d-1}\,.
\]
This implies that the right hand side of (\ref{ibp}) is bounded below by a constant.
Therefore, (\ref{dOX}) yields
\bel{BB1}
\ds{d\over dt} \meas \Big(\Omega^\xi(t)\Big)~\geq ~  \vp'(1)\cdot\meas
\bigl(\Omega^\xi(t)\bigr) - {\omega_{d-1} M d \over d-1}.
\eeq
By Gronwall's inequality, for all times $t\geq 0$  we conclude  that 
\bel{LB1}\meas \Big(\Omega^\xi(t)\Big)~\geq~e^{ \vp'(1)\cdot t}\, \meas(\Omega_0) - {\omega_{d-1}M d \over d-1}\cdot t\,.\eeq
This a priori lower bound on the measure of the set $\Omega^\xi(t)$ shows that
approximate controllability cannot be achieved.
\endproof
\v
{\bf Remark.} Assume that the divergence of the vector field $\bfv(\cdot, \xi)$ remains negative 
for $x$ close to $\xi$, that is
\bel{ndiv}
\vp'(r)+{d-1\over r}\,\vp(r)~\leq~0\qquad\qquad\forall 0<r<\bar r\,,\eeq
for some $\bar r>0$. In this case, the Global Approximate Confinement property implies
\bel{EC2} \limsup_{r\to 0+} ~r^{d-1}\vp(r)~=~+\infty.\eeq
Indeed, one can replace the estimate (\ref{dOX}) with
\bel{dO2}
\bega{rl}\ds{d\over dt} \meas \Big(\Omega^\xi(t)\Big)&\ds
=~\int_{\Omega^\xi(t)\cap B(\xi(t), \bar r)} \div 
\,\bfv\, dx+\int_{\Omega^\xi(t)\setminus B(\xi(t), \bar r)} \div \,\bfv\, dx
\\[4mm]
&\ds \geq~ \int_{\delta<|x|< \bar r} \div 
\,\bfv\, dx+ \vp'(\bar r)\cdot\meas
\bigl(\Omega^\xi(t)\bigr)\\[4mm]
&\ds=~\omega_{d-1} \Big[ \bar r^{d-1} \vp(\bar r) - \delta^{d-1}\vp(\delta)\Big]+\vp'(\bar r)\cdot\meas
\bigl(\Omega^\xi(t)\bigr)\\[4mm]
&\ds\geq ~-\omega_{d-1} \delta^{d-1}\vp(\delta)+\vp'(\bar r)\cdot\meas
\bigl(\Omega^\xi(t)\bigr)
\,.\enda
\eeq
If (\ref{EC2}) fails, then 
\[
{d\over dt} \meas \Big(\Omega^\xi(t)\Big)~\geq~\vp'(\bar r)\cdot\meas
\bigl(\Omega^\xi(t)\bigr)-C_1
\]
for some constant $C_1$, and it 
leads again a priori lower bound on the measure of the set $\Omega^\xi(t)$.

\section{Proof of the  sufficient condition} 
\label{sec:3}
\setcounter{equation}{0}

Aim of this section is to provide a proof of Theorem 2.
As a preliminary, 
consider a bounded open set $\Omega\subset \R^d$ with $\C^2$ 
boundary $\Sigma=\partial \Omega$. On the complement $\R^d\setminus\Sigma$
we consider the vector field
\bel{vx}
\bfv(x)~\doteq~\int_\Sigma\vp(|x-\xi|)\,  {x-\xi\over|x-\xi|}\, d\sigma(\xi),\eeq
where $\sigma$ denotes the $(d-1)$-dimensional surface measure on $\Sigma$.
Since $\Sigma$ has $\C^2$ regularity,  for every $x$ sufficiently close to $\Sigma$ there exists a unique perpendicular 
projection $y_x\in \Sigma$ such that 
\bel{proj}|x-y_x|~=~d(x, \Sigma)~\doteq~\min_{y\in\Sigma} |x-y|\,.\eeq
To fix the ideas, assume that this perpendicular 
 projection $x\mapsto y_x$ is well defined whenever
$d(x, \Sigma)< r_0$, for some curvature radius $r_0>0$.  In the following, 
$\bfn_x= {x-y_x\over |x-y_x|}$ denotes the unit normal 
to the surface $\Sigma$ at the point $y_x$.
\v
{\bf Lemma 1.} {\it Let the function $\vp$ satisfy {\bf (A1)-(A2)}.
Then the vector field $\bfv$ in (\ref{vx}) satisfies}
\bel{vb1}
\lim_{d(x,\Sigma)\to 0} \la \bfn_x\,,\, \bfv(x)\ra ~=~+\infty.
\eeq
\begin{figure}[htbp]
\centering
\includegraphics[scale=0.45]{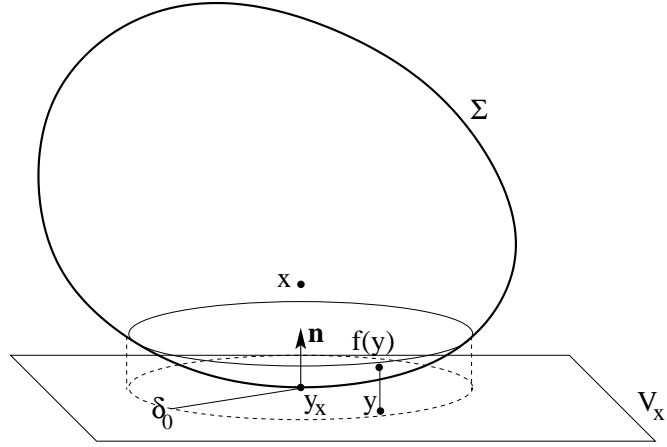}
    \caption{\small  The portion of $\Sigma$ near the point $y_x$ can be represented as
    the graph of a function $f$, as in (\ref{Sig1}).}
\label{f:di10}
\end{figure}
\quad\\
{\bf Proof.} {\bf 1.}
Consider any point $x$ sufficiently close to $\Sigma$ so that the 
perpendicular projection $y_x\doteq \pi(x)$ at
(\ref{proj}) is well defined.   
Let $V_x$ be the hyperplane tangent to  $\Sigma$ at $y_x$ and let 
$\bfn_x$ be the unit normal vector.
As shown in Fig.~\ref{f:di10}, 
in a neighborhood of $y_x$, the surface $\Sigma$ can be expressed as the 
graph of a function $f:V_x\mapsto\R$. More precisely, call $\ve= d(x,\Sigma)=|x-y_x|$.
Without loss of generality, we can choose a system of coordinates 
such that  $y_x=0$.
Notice that, by the regularity and compactness of the surface $\Sigma$, we can assume that 
the radius $\delta_0$ of the ball where the function $f$ is defined is independent of 
$y_x\in \Sigma$. Moreover, the $\C^2$ norm of $f$ remains uniformly bounded.
By construction we have
\bel{fs}
f(0)~=~0, \qquad \nabla f(0)~=~0,\qquad \|f\|_{\C^2(B(0, \delta_0))}~\leq~C_0,\eeq
for some uniform constant $C_0$.
Defining the constant $\kappa = \sqrt {2/C_0}>0$, by (\ref{fs}) one has
the implication
\bel{kp}
|y-y_x|~<~\kappa |x-y_x|^{1\over2}\qquad\implies\qquad \la \bfn_x\,,~
x-y-f(y)\bfn_x\ra~\geq~0.\eeq

\begin{figure}[htbp]
\centering
\includegraphics[scale=0.3]{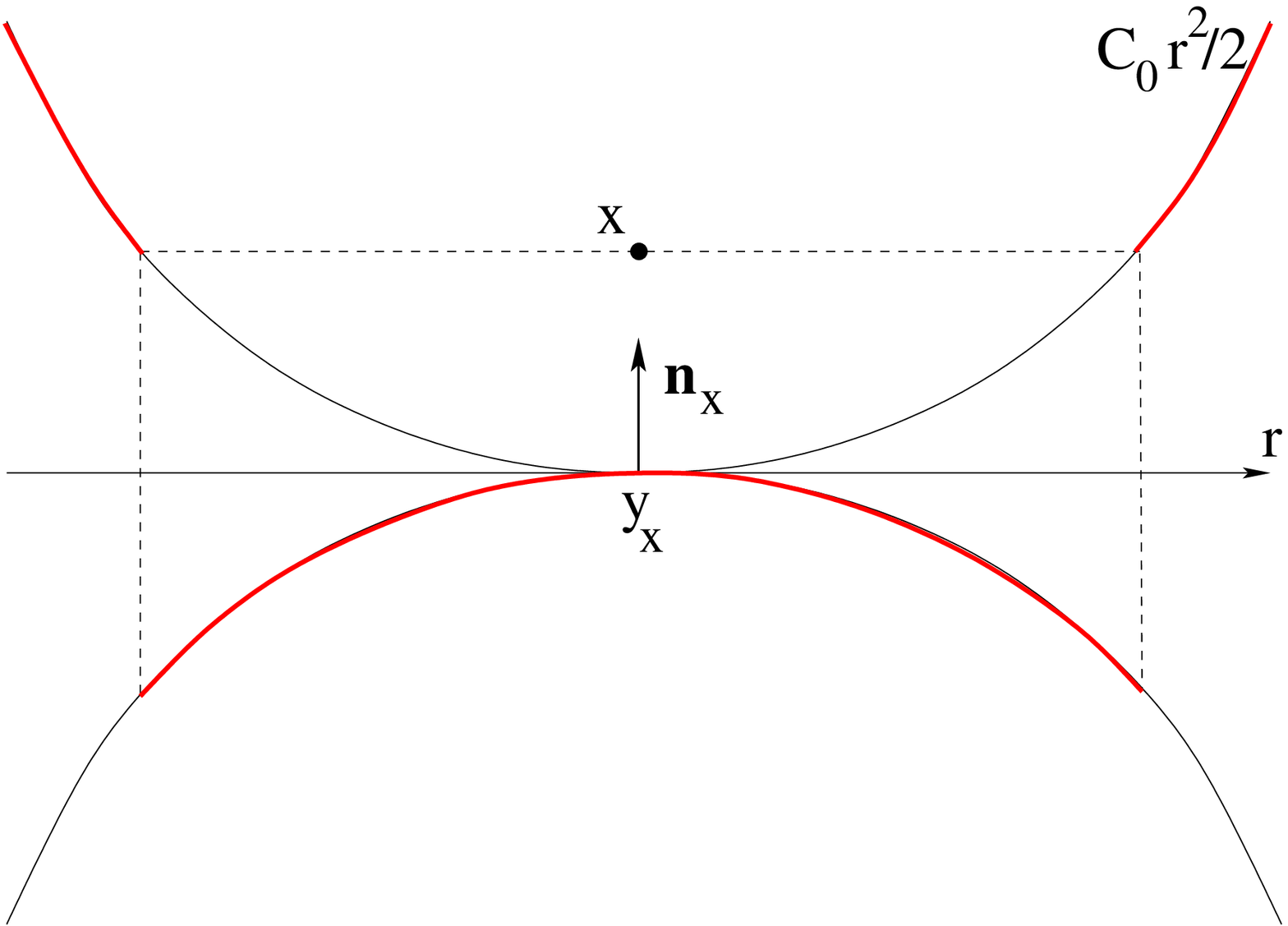}
    \caption{\small  The estimates (\ref{nv1})-(\ref{nv2}).}
\label{f:di12}
\end{figure}

\v
{\bf 2.} Next,
consider the decomposition $\Sigma=\Sigma_1\cup\Sigma_2\cup\Sigma_3$,
where
\bel{S123}\left\{
\bega{rl}\Sigma_1&=~\bigl\{y+f(y)\bfn_x\,;~~|y|\leq\kappa |x|^{1/2}\bigr\},\\[3mm]
\Sigma_2&=~\bigl\{y+f(y)\bfn_x\,;~~\kappa |x|^{1/2}<y<\delta_0\bigr\},\\[3mm]
\Sigma_3&=~\Sigma\setminus(\Sigma_1\cup\Sigma_2).\enda\right.\eeq
Based on (\ref{S123}), we 
shall estimate the vector field $\bfv$ by splitting the integral in (\ref{vx}) in three
parts:
\bel{v123}
\bfv_i(x)~=~\int_{\Sigma_i} \vp(|x-\xi|) \, {x-\xi\over |x-\xi|}\, d\sigma(\xi),\qquad\qquad
i=1,2,3,\eeq
so that
$$\bfv(x) ~=~\bfv_1(x)+\bfv_2(x)+\bfv_3(x).$$ 
In the following, for $y\in V_x$, we use the bound 
$\bigl| f(y)\bigr|~\leq~{C_0\over 2}|y|^2$
and the identities
\bel{xyid}|x|\,=\,\ve\,,\qquad  |y|\,=\,r\,,\qquad |x-y|\,=\, \sqrt{\ve^2+r^2}\,,
\qquad \bfn_x\,=\,{x\over
|x|}\,.\eeq
As long as  $|y|\leq \delta_0$, the above implies
\bel{ulb}
|x-y-f(y)\bfn_x|^2~\in~\Big[ \ve^2 + (1-C_0\ve)r^2\,,~\ve^2 + (1+C_0\ve)r^2 + C_0^2 r^4/4\Big].
\eeq
Using (\ref{ulb}) and the monotonicity of $\vp$ we obtain the estimates
\bel{nv1}
\bega{rl}
 \la \bfn_x\,,\, \bfv_1(x)\ra &\geq~\ds c_0 \int_0^{\kappa \ve^{1/2}}
r^{d-2}\Big(\ve-{C_0\over2}r^2\Big)\cdot {\vp\Big(\sqrt{ (\ve +C_0 r^2/2)^2+r^2}\Big)\over 
\sqrt{ (\ve +C_0 r^2/2)^2+r^2}}\, dr\\[4mm]
&\geq~ \ds c_0 \int_0^{\kappa \ve^{1/2}}
r^{d-2}\Big(\ve-{C_0\over2}r^2\Big)\cdot{\vp\Big(\sqrt{ \ve^2 +(1+2C_0\ve)r^2}\Big)\over 
\sqrt{\ve^2 +(1+2C_0\ve)r^2}}\, dr\\[4mm]
&\geq~\ds {3\over 4}c_0\int_0^{\kappa \ve^{1/2}\over 2}
\ve r^{d-2}\cdot{\vp\Big(\sqrt{ \ve^2 +(1+2C_0\ve)r^2}\Big)\over 
\sqrt{\ve^2 +(1+2C_0\ve)r^2}}\, dr\,,
\enda
\eeq
\bel{nv2}
\bega{rl}
\left| \la \bfn_x\,,\, \bfv_2(x)\ra \right|&\leq~\ds c_0 \int_{\kappa \ve^{1/2}}^{\delta_0}r^{d-2}
\vp\Big(\sqrt{ \ve^2 +(1-C_0 \ve)r^2}\Big)\cdot {\ve+C_0 r^2/2\over 
\sqrt{ \ve^2 +(1-C_0 \ve)r^2}}\, dr\\[4mm]
&\leq~\ds C_1\cdot \int_{\kappa \ve^{1/2}}^{\delta_0} 
{\vp\Big(\sqrt{ \ve^2 +(1-C_0 \ve)r^2}\Big)\over\sqrt{ \ve^2 +(1-C_0 \ve)r^2} }\cdot r^{d}
\, dr\,,
\enda
\eeq
\bel{nv3}
\left| \la \bfn_x\,,\, \bfv_3(x)\ra\right| ~\leq~  C\,,\eeq
for some constants $C,C_1, c_0>0$. 
 Performing the change of variable $s=r\,\sqrt{1+2C_0\ve}$ in (\ref{nv1}), one finds
\[
\bega{rl}
 \la \bfn_x\,,\, \bfv_1(x)\ra&\geq~\ds C_2\cdot\int_{0}^{{1\over 2} \kappa \ve^{1/ 2}}\ve s^{d-2}\cdot {\vp\left(\sqrt{\ve^2+s^2}\right)\over \sqrt{\ve^2+s^2}}~ds\\[4mm]
 &\geq~\ds{4C_2\over \kappa^2}\cdot \int_{0}^{{1\over 2}\kappa \ve^{1/2}}s^{d}\cdot {\vp\left(\sqrt{\ve^2+s^2}\right)\over \sqrt{\ve^2+s^2}}~ds
 \enda
\]
for some constant $C_2>0$. Setting $t=\sqrt{\ve^2+s^2}$, we estimate 
\[
\bega{rl}
 \la \bfn_x\,,\, \bfv_1(x)\ra&\geq~\ds {4C_2\over \kappa^2}\cdot\int_{\ve}^{\sqrt{\ve^2+{1\over 4}\kappa^2\ve}} \left(t^2-\ve^2\right)^{d-1\over 2}\cdot\vp(t)~dt\\[4mm]
 &\geq~\ds {4C_2\over \kappa^2}\cdot\int_{2\ve}^{{1\over 2} \kappa \ve^{1/ 2}} \left(t^2-\ve^2\right)^{d-1\over 2}\cdot\vp(t)~dt\\[4mm]
 &\geq~\ds {4C_2\over \kappa^2}\cdot\left(3\over 4\right)^{d-1\over 2}\cdot \int_{2\ve}^{{1\over 2} \kappa \ve^{1/ 2}}~ t^{d-1}\cdot\vp(t)~dt\,.
 \enda
\]
Similarly, using the variable $s=\sqrt{\ve^2+(1-C_0\ve)r^2}$ in (\ref{nv2}), we have 
\begin{eqnarray*}
\int_{\kappa \ve^{1/2}}^{\delta_0} 
{\vp\Big(\sqrt{ \ve^2 +(1-C_0 \ve)r^2}\Big)\over\sqrt{ \ve^2 +(1-C_0 \ve)r^2} }\cdot r^{d}
\, dr&=&{1\over 1-C_0\ve}\cdot\ds\int^{\sqrt{\delta^2+\ve(\ve-C_0\delta^2})}_{\sqrt{\ve^2+(1-C_0\ve)\kappa^2 \ve}}\varphi(s)\cdot r^{d-1}~ds\\[4mm]
&\leq& {1\over (1-C_0\ve)^{d+1\over 2}}\cdot \ds\int^{\sqrt{\delta^2+\ve(\ve-C_0\delta^2})}_{\sqrt{\ve^2+(1-C_0\ve)\kappa^2 \ve}}\varphi(s)\cdot s^{d-1}~ds\,.
\end{eqnarray*}
Thus, for $\ve>0$  sufficiently small, it holds
\[
\bega{rl}
\left| \la \bfn_x\,,\, \bfv_2(x)\ra \right|~\leq~\ds C_3\cdot\int^{\delta_0}_{{1\over 2}\kappa \ve^{1/ 2}}~s^{d-1}\cdot\vp(s)~ds
\enda
\]
for some constant $C_3>0$. Setting $\tilde{\ve}=2\ve$ we obtain 
\bel{nv1-e2}
 \la \bfn_x\,,\, \bfv_1(x)\ra~\geq~\ds C_4 \cdot \int_{\tilde{\ve}}^{\kappa_1 \tilde{\ve}^{1/ 2}}~s^{d-1}\cdot\vp(s)~ds
\eeq
and
\bel{nv2-e2}
\left| \la \bfn_x\,,\, \bfv_2(x)\ra \right|~\leq~C_3\cdot \int^{\delta_0}_{\kappa_1 \tilde{\ve}^{1/ 2}}~s^{d-1}\cdot\vp(s)~ds\,,
\eeq
for $\kappa_1={\kappa\over 2\sqrt{2}}$ and some constant $C_4>0$. In particular, 
\[
{ \la \bfn_x\,,\, \bfv_1(x)\ra\over \left| \la \bfn_x\,,\, \bfv_2(x)\ra \right|}~\geq~{C_4\over C_3}\cdot{\ds\int_{\tilde{\ve}}^{\kappa_1 \tilde{\ve}^{1/ 2}}~s^{d-1}\cdot\vp(s)~ds\over \ds\int^{\delta_0}_{\kappa_1 \tilde{\ve}^{1/ 2}}~s^{d-1}\cdot\vp(s)~ds} ~=~{C_4\over C_3}\cdot \left[ {\ds \int_{\tilde{\ve}}^{\delta_0}~s^{d-1}\cdot\vp(s)~ds\over  \ds \int^{\delta_0}_{\kappa_1 \tilde{\ve}^{1/ 2}}~s^{d-1}\cdot\vp(s)}-1\right]\,.
\]
If the assumption {\bf(A2)} holds, then one has
\[
\limsup_{r\to 0+}~r^{d}\cdot \vp(r)~=~+\infty\,.
\]
This implies 
\[
\lim_{r\to 0+}\int_{r}^{\delta_0}s^{d-1}\cdot \vp(s)~ds~\geq~\limsup_{r\to 0+}\int_{r}^{2r}s^{d-1}\cdot \vp(s)~ds~\geq~{2^d-1\over d2^d}\cdot \limsup_{r\to 0+}~ (2r)^{d}\cdot \vp(2r)~=~+\infty\,.
\]
Applying L'Hopital's rule and the assumption {\bf(A2)}, we obtain
\begin{eqnarray*}
\lim_{\tilde{\ve}\to 0+}~{\ds\int^{\delta_0}_{\kappa_1 \tilde{\ve}^{1/ 2}}~s^{d-1}\cdot\vp(s)~ds\over\ds\int_{\tilde{\ve}}^{\delta_0}~s^{d-1}\cdot\vp(s)~ds }&=&{\kappa^d_1\over 2}\cdot \lim_{\tilde{\ve}\to 0+}~{\tilde{\ve}^{d/2-1}\cdot\vp(\kappa_1\tilde{\ve}^{1/2})\over \tilde{\ve}^{d-1}\vp(\tilde{\ve})}\\[4mm]
&=&{\kappa^d_1\over 2}\cdot \lim_{\tilde{\ve}\to 0+}~{\tilde{\ve}^{d/2}\cdot\vp(\kappa_1\tilde{\ve}^{1/2})\over \tilde{\ve}^{d}\vp(\tilde{\ve})}~=~0.
\end{eqnarray*}
This yields
\[
\left|{\la \bfn_x\,,\, \bfv_2(x)\ra \over \la \bfn_x\,,\, \bfv_1(x)\ra}\right|~\leq~{C_3\over C_4}\cdot {\ds\int^{\delta_0}_{\kappa_1 \tilde{\ve}^{1/ 2}}~s^{d-1}\cdot\vp(s)~ds\over\ds\int_{\tilde{\ve}}^{\kappa_1 \tilde{\ve}^{1/ 2}}~s^{d-1}\cdot\vp(s)~ds }~\to~0\qquad\mathrm{as}\qquad\ve~\to~0\,.
\]
Finally, observing that 
$$| \la \bfn_x\,,\, \bfv_3(x)\ra|~\leq~
C,$$
the limit behavior (\ref{vb1}) is clear. 
This achieves the proof, because the constants $C_0, \delta_0, \kappa$
are independent of the  point $y_x\in \Sigma$.
\endproof

In the following, 
for $x\notin \Sigma(t)$ we denote by $\pi(t,x)$ the perpendicular projection 
of $x$ on $\Sigma(t)$, and call
\bel{ppro}\bfn(t,x)~\doteq~{x-\pi(t,x)\over |x-\pi(t,x)|}\eeq
the unit normal vector.
\v
{\bf Corollary 1.} {\it Consider a family of compact $\C^2$ surfaces $\Sigma(t)$, 
continuously depending on $t\in [0,T]$, with uniformly bounded curvature.
Define the vector fields 
\bel{bvdef}\bfv(t,x)~\doteq~\int_{\Sigma(t)}\vp(|x-\xi|)\,  {x-\xi\over|x-\xi|}\, d\sigma(\xi).\eeq
Then for any $N$ there exists $\delta>0$ such that, for all $t\in [0,T]$,
\bel{vb5}
d(x,\Sigma(t))<\delta\qquad\implies\qquad  \la \bfn(t,x)\,,\, \bfv(t,x)\ra ~\geq ~N.
\eeq}
\v
Indeed, this follows from Lemma~1, observing that the limit in (\ref{vb1}) is uniform over
all surfaces $\Sigma(t)$.

\v
{\bf Proof of Theorem 2.}

{\bf 1.} 
Let the compact sets $\Omega_1\subset int\,\Omega_0$ be given, with $\Omega_1$ simply connected.  To fix the ideas, assume 
\bel{O10}B(\Omega_1, \rho)~\subset~\Omega_0\,,\eeq
for some radius $\rho>0$.
Given $T>0$ and $0<\ve<\rho$, we can find a decreasing 
family of  compact sets $t\mapsto V(t)$ with 
$\C^2$ boundary, as in (\ref{Vdef})-(\ref{ppp}),  such that
\bel{Vt2}
B\bigl(\Omega_0, \ve\bigr)~\subset~V(0),\qquad\qquad    
B(\Omega_1, \ve/2)~\subset~V(T)~\subset~ B(\Omega_1, 3\ve/4).\eeq
Call  $\Sigma(t)\doteq \partial V(t)$ the boundaries of these sets and define
the vector fields 
\bel{nvt}\bfw(t,x)~\doteq~\delta_0\cdot \int_{\Sigma(t)}\vp(|x-\xi|)\,  {x-\xi\over|x-\xi|}\, d\sigma(\xi).\eeq
Here the constant  $\delta_0>0$ is chosen small enough so that
\bel{dvsm}\delta_0\cdot \int_{\Sigma(t)}\, d\sigma~\leq~1\qquad\qquad\forall t\in [0,T],\eeq
\bel{vsm}
\bigl|\bfw(t,x)\bigr|~<~{\ve\over 8T}\qquad\forall x\in B(\Omega_1, \ve/2), ~~~t\in [0,T].\eeq
\v
{\bf 2.}  
For any point $x_0\in \Omega_0$, denote by $t\mapsto x(t, x_0)$ the solution of
\bel{avode}\dot x(t)~=~\bfw(t, x(t)),\qquad\qquad x(0)=x_0\,.\eeq
We claim that 
\bel{xin}
x(t, x_0)~\in~int \,V(t)\qquad\qquad\forall t\in [0,T].\eeq
To prove (\ref{xin}), let $L$ be
 a Lipschitz 
constant  for the multifunction $t\mapsto\Sigma(t)$, so  that 
the Hausdorff distance between the two boundaries satisfies
\bel{Dst} d_H\Big( \Sigma(s),\, \Sigma(t)\Big)~\leq~L\,(t-s)\qquad \hbox{ for any}\quad 0<s<t<T.
\eeq
For any trajectory $t\mapsto x(t)$ of (\ref{avode}), (\ref{nvt}),   consider the distance
\[
d(t)~=~\hbox{dist}\,\Big( x(t)\,,~\Sigma(t)\Big)\]
 of $x(t)$ from the boundary  $\Sigma(t)=\partial V(t)$. 
By Corollary~1 there is a constant $\ve_1\in \,]0,\ve]$  such that,
for any $x\in int \,V(t)$ with dist$\bigl(x, \Sigma(t)\bigr)\leq \ve_1$, one has 
\bel{paway}
\left\langle \bfn(t,x),\int_{\Sigma(t)} \vp(|x-\xi|)\,{x-\xi\over 
|x-\xi|}\, d\sigma(\xi)\right\rangle~>~{L\over \delta_0}\,.
\eeq
In view of (\ref{paway}) and (\ref{Dst}), if  $d(x(t),\Sigma(t))\leq \ve_1$,  then the time derivative 
of the distance $d(t)$ satisfies
\bel{ddist}\bega{l}
\dot{d}(t)~
\geq~\ds-L+\Big\langle \bfn(t,x(t)) ,\, \bfw(t, x(t))\Big\rangle\\[3mm]
\qquad \ds =~-L+ 
\left\langle \bfn(t, x(t)), \,\delta_0\cdot\int_{\Sigma(t)} \vp(|x(t)-\xi|)\,{x(t)-\xi\over 
|x(t)-\xi|}\, d\sigma(\xi)\right\rangle\\[3mm]
\qquad \ds >~-L + \delta_0 \cdot {L\over \delta_0}~=~0\,.\enda
\eeq
If $x(t)$ is a trajectory starting inside $\Omega_0$, then $d(0)\geq\ve \geq \ve_1$.
By (\ref{ddist}) we thus have $d(t)\geq \ve_1$ for all $t\in [0,T]$.
\v
We conclude this step by observing that,  for every $x_0\in \ov B(\Omega_1, \ve/4)$, 
by (\ref{vsm}) the corresponding trajectory satisfies
\bel{afix}
|x(t,x_0)- x_0|~\leq~{\ve t\over 8T}\qquad\qquad \forall t\in [0,T].
\eeq
\quad\\
{\bf 3.} 
Relying on  the approximation procedure developed in
\cite{BZ}, we claim that there exists a Lipschitz control $t\mapsto \xi(t)$ for (\ref{G})-(\ref{ODE1})
that
produces almost the same trajectories as (\ref{avode}).   
More precisely,
calling 
$t\mapsto x^\xi(t, x_0)$ the solution to
\bel{xv2}\dot x~=~\bfv(x,\xi(t)),\qquad\qquad x(0)=x_0\,,\eeq
for every  $t\in [0,T]$ and $x_0 \in \Omega_0$ one has
\bel{close}
\bigl| x^\xi(t,x_0) - x(t, x_0)\big|~<~{\ve\over 8}\,.\eeq
Toward this goal, for any $t\in [0,T]$, define $\mu^t$ to be the $(d-1)$-dimensional measure 
supported on $\Sigma(t)$,
so that
$$\mu^t(A)~=~\int_{A\cap\Sigma(t)} \, d\sigma$$
 for every open set $A\subset\R^d$.  
By (\ref{dvsm}) it follows
 $$\delta_0\,\mu^t(\R^d)~=~\delta_0
 \cdot[\hbox{surface area of}~\Sigma(t)]~\leq~1\qquad\qquad \forall t\in [0,T].$$
 We now choose a point  $\bar x\in \R^d$ very far from the origin and define the probability measure
 $$\tilde\mu^t~=~\delta_0\mu^t + \Big( 1-\delta_0\,\mu^t(\R^d)\Big)m_{\bar x}\,$$
 where $m_{\bar x}$ denotes a unit Dirac mass at $\bar x$.

 Notice that, as $|\bar x|\to +\infty$, by the second limit in (\ref{vplim}) the  vector
 \bel{TG}
\Tilde \bfw(t,x)~\doteq~\int
\vp(|x-\xi|)\,{x-\xi\over |x-\xi|}\,d\tilde\mu^t(\xi)\eeq
 approaches $\delta_0 \bfw(t,x)$, uniformly on compact subsets of $\R^d\setminus\Sigma(t)$.

Given an integer $n\geq 1$, we split the interval $[0,T]$ into $n$ equal subintervals,
inserting the points 
$t_i~=~iT/n$, $i=0,1,\ldots,n$.
For each $i$,
the probability measure $\tilde \mu^{t_i}$  can now be approximated by 
the sum of $N$ equal masses, say located at $\xi_{i1},\ldots, \xi_{iN}$.
Defining the time step  $h \doteq {T\over n\,N}$, we then consider the control function
\bel{xapp}\xi(t)~=~\xi_{ij}\qquad \hbox{if}\quad t_i +(j-1)h< t ~\leq~t_i+j h\,.\eeq
The same arguments used in \cite{BZ} now show that, as $n,N\to \infty$, by suitably choosing the points 
$\xi_{ij}$,
 trajectories of the ODE (\ref{ODE1}), (\ref{xapp}) 
converge to the corresponding trajectories of $\dot x = \tilde \bfv(t,x)$. Moreover, the convergence 
is  uniform
for all initial data in the compact set $\Omega_0\subset \,\R^d\setminus\Sigma(0)$.

Finally, we can replace the piecewise constant function $\xi(\cdot)$ by a Lipschitz function $\Tilde\xi(\cdot)$.
If $\|\Tilde\xi-\xi\|_{\L^1}$ is sufficiently small, the corresponding trajectories still satisfy the 
same estimate (\ref{close}).
\v
{\bf 4.} 
Recalling that 
$$x(T, x_0) ~\in~V(T)~\subseteq~B(\Omega_1, \, 3\ve/4),$$
by (\ref{close}) we now conclude
\bel{Oxi}\Omega^\xi(T)~=~\bigl\{  x^\xi(T,x_0)\,;~~x_0\in \Omega_0\bigr\}
~\subseteq~B(\Omega_1, \, \ve).\eeq
This establishes the second inclusion in (\ref{DHT}).

To prove the first inclusion, consider the continuous map $x_0\mapsto x^\xi(T, x_0)$
from the compact set $\ov B(\Omega_1,\, \ve/4)$ into $\R^d$.
By (\ref{afix}) and (\ref{close}) it follows
\bel{cl1}\bigl|x^\xi(T, x_0)-x_0\bigr|~\leq~{\ve\over 4}\qquad\forall x_0\in \ov B\Big(\Omega_1,\, 
{\ve\over 4}\Big).
\eeq
For any  given $y\in\Omega_1$, define the continuous map 
$g^y:\overline{B}(0,\varepsilon/4)\mapsto\R^d$  by setting
\[
g^y(z)~\doteq~x^{\xi}(T,y-z)-(y-z)\qquad\forall z\in \overline{B}(0,\ve/4)\,.
\]
By (\ref{cl1}) one has
 $$g^y(\overline{B}(0,\varepsilon/4))~\subseteq~\overline{B}(0,\varepsilon/4).$$ 
Therefore, Brouwer's fixed point theorem implies
\[
g^y(z_0)~=~z_0\qquad\mathrm{for~some}~z_0\in \overline{B}(0,\ve/4).
\]
This yields 
\[
y~=~x^\xi(T, x_0)\qquad\mathrm{with}\qquad x_0~=~y-z_0\in\ov B(y, \ve/4).
\]
Hence $\Omega_1\subseteq \Omega^\xi(T)$.
\v
{\bf 5.} Finally,
to pass from approximate controllability to exact controllability one can 
split the interval $[0,T]$, inserting an increasing sequence of times $\tau_j$ with $\tau_j\to T-$ as $j\to\infty$.   Then construct Lipschitz controls $t\mapsto \xi(t)$ 
on each subinterval $[\tau_{j-1}, \tau_j]$   such that the corresponding sets
$\Omega^\xi(\tau_j)$ satisfy
$$\ov B(\Omega_1, 2^{-j})~\subseteq~\Omega^\xi(\tau_j)~\subseteq~\
\ov B(\Omega_1, 2^{1-j}).$$
\endproof




\section{Approximating a sweeping process} 
\label{sec:4}
\setcounter{equation}{0}  

The key tool for the proof of Theorem 3 is the following lemma, which 
improves on Lemma 1 under the stronger assumption {\bf (A2$'$)}.
\v
{\bf Lemma 2.} {\it 
Let $\Omega\subset \R^d$ be a compact set with $\C^2$ 
boundary $\Sigma=\partial \Omega$. 
Let $\bfv$ be the vector field in (\ref{vx}).  
If the function $\vp$ satisfies {\bf (A1)-(A2$'$)}, then
\bel{vbig}
\bigl|\bfv(x)\bigr|~\to~+\infty\qquad\qquad \hbox{as} \quad d(x, \Sigma)~\to~0,\eeq
\bel{vperp}\left| {\bfv(x)\over \bigl|\bfv(x)\bigr|}- {x-\pi( x)\over \bigl|x-\pi( x)\bigr|}\right|~\to~0
\qquad\hbox{as} \quad d(x, \Sigma)~\to~0.\eeq
}
\v
{\bf Proof.} {\bf 1.}
Consider any point $x$ sufficiently close to $\Sigma$ so that the 
perpendicular projection $y_x\doteq \pi(x)$ at
(\ref{proj}) is well defined.   As in the proof of Lemma~1,
in a neighborhood of $y_x$, the surface $\Sigma$ can be expressed as the 
graph of a function $f:V_x\mapsto\R$. 
More precisely, call $\ve= d(x,\Sigma)=|x-y_x|$.
Then, given ${1\over2}<\alpha<1$, we can write $\Sigma=\Sigma_1\cup\Sigma_2\cup\Sigma_3$, where
\bel{Sig1}\left\{\bega{rl}
\Sigma_1&=~\Big\{ y+f(y)\, \bfn_x\,;~~y\in V_x\,,~|y|\leq \ve^\alpha\Big\},\\[3mm]
\Sigma_2&=~\Big\{ y+f(y)\, \bfn_x\,;~~y\in V_x\,,~\ve^\alpha<|y|\leq \delta_0\Big\},
\\[3mm]
\Sigma_3&=~\Sigma\setminus(\Sigma_1\cup\Sigma_2)\,.
\enda\right. \eeq
Notice that, by the regularity and compactness of the surface $\Sigma$, we can assume that 
the radius $\delta_0$ of the ball where the function $f$ is defined is independent of 
$y_x\in \Sigma$. Moreover, the $\C^2$ norm of $f$ remains uniformly bounded.

Without loss of generality, in the following computations  we 
shall assume $y_x=0\in \R^d$.
By construction we again have the bounds (\ref{fs}), valid
for some constant $C_0$, uniform w.r.t.~$y_x\in \Sigma$.
We shall estimate the vector field 
$$\bfv(x) = \bfv_1(x)+\bfv_2(x)+\bfv_3(x)$$ 
by splitting the integral (\ref{vx}) in three
parts, as in (\ref{v123}). Notice, however, that now we refer to the different decomposition
(\ref{Sig1}) of the surface $\Sigma$. 
\v
{\bf 2.} Calling $J(y)=\sqrt{1+|\nabla f(y)|^2}$ the Jacobian determinant of the map $y\mapsto y+f(y)\bfn_x$ from $V_x\cap B(y_x,\delta_0)$ into $\Sigma$, 
we have 
\bel{v1a}
\bfv_1(x)~=~
\int_{|y|<\ve^\alpha}\vp\bigl(|x-y-f(y) \bfn_x|\bigr)\cdot
{x-y-f(y) \bfn_x\over|x-y-f(y) \bfn_x|}\, J(y)\, dy.
\eeq
We write
$$\bfv_1(x)~=~\bfv_{11}(x)+\bfv_{12}(x),$$
where
\bel{v11}\bfv_{11}(x)~=~
\int_{|y|<\ve^\alpha}\vp\bigl(|x-y|\bigr)\,
{x-y\over|x-y|}\,dy\,,\qquad\qquad \bfv_{12}(x)~=~ \bfv_1(x)-\bfv_{11}(x).
\eeq
Notice that $\bfv_{11}(x)$  is  a vector parallel to $\bfn_x$ and is computed as
\bel{v11e}
\bega{rl}
\bfv_{11}(x)~=~ \left(c_0\ds\int_0^{\ve^\alpha}r^{d-2} \vp\Big(\sqrt{\ve^2+r^2}
\Big)\cdot {\ve\over \sqrt{\ve^2+r^2}}\, dr\right)\cdot{\bfn _x} \,,
\enda
\eeq
for some constant $c_0>0$\,. 
Hence, in order to obtain the limit (\ref{vperp}), let's first prove that
\bel{lim03}
{\bigl|{\bf v}_{12}(x)\bigr|
\over \big|{\bf v}_{11}(x)\big|}~\to~0\qquad\hbox{as}\quad \ve\to 0\,.\eeq

The vector $\bfv_{12}(x)$ satisfies
 \bel{v12}\bega{rl}
 \bfv_{12}(x)&=~\ds\int_{|y|<\ve^\alpha}\Big[ \vp\bigl(|x-y-f(y) \bfn_x|\bigr)-
 \vp\bigl(|x-y|\bigr)\Big] \cdot
{x-y-f(y) \bfn_x\over|x-y-f(y) \bfn_x|}\, J(y)\, dy\\[4mm]
 &\qquad \ds+\int_{|y|<\ve^\alpha}\vp\bigl(|x-y|\bigr)\cdot
 \left\{{x-y-f(y) \bfn_x\over|x-y-f(y) \bfn_x|}-
{x-y\over|x-y|}\right\} \, J(y)\, dy
\\[4mm]
 &\qquad \ds+
 \int_{|y|<\ve^\alpha}\vp\bigl(|x-y|\bigr)\cdot
{x-y\over|x-y|}\, \bigl[J(y)-1\bigr]\, dy
\\[4mm]
 &=~A_1+A_2+A_3
 \,.\enda
\eeq


In the following, recalling (\ref{fs}), we use
the bounds
\bel{BBB}\bigl| f(y)\bigr|~\leq~{C_0\over 2}|y|^2,\qquad\quad
|J(y) - 1|~=~\O(1)\cdot |y|^2,\eeq
and the identities (\ref{xyid}).
Since the function $\vp$ is decreasing and $r\leq\ve^\alpha$, 
using (\ref{ulb}) one obtains  the estimate
\bel{A1}
\bega{rl}|A_1|&\ds\leq~C\cdot \int_0^{\ve^\alpha} r^{d-2} \, \left[\vp\Big(\sqrt{\ve^2+(1-C_0\ve)r^2}
\Big)-\vp\Big(\sqrt{\ve^2+r^2}
\Big)\right] \, dr\,\\[4mm]
&\ds\qquad\qquad\qquad\qquad+
C\cdot \int_0^{\ve^\alpha} r^{d-2} \, \left[\vp\Big(\sqrt{\ve^2+r^2}
\Big)-\vp\Big(\sqrt{\ve^2+(1+2C_0\ve)r^2}
\Big)\right]\,dr\\[4mm]
&\doteq~ A_{11}+A_{12}\,,
\enda
\eeq
for some constant $C$ and all $\ve>0$ sufficiently small.
In addition, we have
\bel{A1__}
\bega{rl} \ds\left|{x-y-f(y) \bfn_x\over|x-y-f(y) \bfn_x|}-
{x-y\over|x-y|}\right|&\leq~\ds\left|{|x-y|-|x-y-f(y) \bfn_x|\over|x-y|}\right|+{|f(y)|\over|x-y|}\\[4mm]
&\leq~2\cdot\ds{|f(y)|\over|x-y|}~\leq~ C_0\cdot\ds{|y|^2\over|x-y|}\,.
\enda
\eeq
Consequently,
\bel{A2}
\bega{rl}
|A_2|&
\leq~C\cdot \ds\int_{|y|<\ve^\alpha} \vp\left( {|x-y|}\right)\cdot {|y|^2\over |x-y|}\, dy~=~C\cdot c_0\ds\int_0^{\ve^\alpha} r^{d-2} \vp\Big(\sqrt{\ve^2+r^2}
\Big) \cdot {r^2\over \sqrt{\ve^2+r^2}}\, dr\\[4mm]
&\leq~C\cdot c_0\ds\int_0^{\ve^\alpha} r^{d-2} \vp\Big(\sqrt{\ve^2+r^2}
\Big) \cdot {\ve^{2\alpha}\over \sqrt{\ve^2+r^2}}\, dr
\enda
\eeq
and 
\bel{A3}
\bega{rl}|A_3|&\leq~C\cdot c_0\ds\int_0^{\ve^\alpha}r^{d-2}\vp\Big(\sqrt{\ve^2+r^2}
\Big) \, r^2\, dr~\leq~C\cdot c_0\ds\int_0^{\ve^\alpha} r^{d-2} \vp\Big(\sqrt{\ve^2+r^2}
\Big) \cdot {\ve^{2\alpha}\over \sqrt{\ve^2+r^2}}\, dr
\enda
\eeq
for a suitable constant $C$.

Since we are choosing $\alpha>1/2$, comparing (\ref{A2}) and (\ref{A3}) with (\ref{v11e}),
it is clear that
\bel{Asm}
{|A_2|+|A_3|\over |\bfv_{11}(x)|}~\leq~2C\ve^{2\alpha-1}~\to ~0\qquad \hbox{as}\quad \ve\to 0.\eeq

Proving a similar estimate for $A_1$ requires more work. 
Performing the variable change
\[
s~=~\sqrt{1-C_0\ve}\cdot r\,,
\]
one obtains
\bel{A1-1}
\bega{rl}\ds 0&\leq~\ds\int_0^{\ve^\alpha} r^{d-2} \, \left[\vp\Big(\sqrt{\ve^2+(1-C_0\ve)r^2}
\Big)-\vp\Big(\sqrt{\ve^2+r^2}
\Big)\right] \, dr \\[4mm]
\quad&\ds\leq~(1-C_0\ve)^{-{d-1\over 2}}\int_{0}^{\ve^{\alpha}\cdot \sqrt{1-C_0\ve}}s^{d-2}\vp\left(\sqrt{\ve^2+s^2}\right)~ds-\int_{0}^{\ve^{\alpha}}r^{d-2}\vp\Big(\sqrt{\ve^2+r^2}
\Big)~dr\\[4mm]
&\ds\leq~\left[(1-C_0\ve)^{-{d-1\over 2}}-1\right]\,\int_{0}^{\ve^\alpha}r^{d-2}\vp\left(\sqrt{\ve^2+r^2}\right)~dr
\\[5mm]
&\ds\leq~C_1\ve\,\int_{0}^{\ve^\alpha}r^{d-2}\vp\left(\sqrt{\ve^2+r^2}\right)~dr ,
\enda
\eeq
for some constant $C_1$. 
Recalling (\ref{A1}) and comparing with (\ref{v11e}), we thus obtain
\bel{A11}\bega{rl}
A_{11}&\leq~\ds C_2\ve\int_{0}^{\ve^\alpha}r^{d-2}\vp\left(\sqrt{\ve^2+r^2}\right)  \,dr \\[4mm]
&\leq~\ds  C_2 \sqrt{\ve^2 + \ve^{2\alpha}}\cdot \int_{0}^{\ve^\alpha}r^{d-2}\vp\left(\sqrt{\ve^2+r^2}\right) \,{\ve\over \sqrt{\ve^2+r^2}} \,dr\\[4mm]
&=~\ds {C_2\over c_0} \sqrt{\ve^2 + \ve^{2\alpha}}\cdot |\bfv_{11}(x)|~\leq~{2C_2\ve^{\alpha}\over c_0}\cdot \big|{\bf v}_{11}(x)\big|\enda
\eeq
for some constant $C_2$.  A similar argument yields
\bel{A1-2}
\bega{rl}
A_{12}&\leq~\ds C\cdot\left[\int_{0}^{\ve^{\alpha}}r^{d-2}\vp\left(\sqrt{\ve^2+r^2}\right)~dr~\ds-\left({1\over 1+2C_0\ve}\right)^{d-1\over 2}\cdot\int_{0}^{\ve^{\alpha}\sqrt{1+2C_0\ve}}r^{d-2}\vp\left(\sqrt{\ve^2+r^2}\right)~dr\right]\\[4mm]
&\leq~\ds C_3\ve \int_{0}^{\ve^{\alpha}}r^{d-2}\vp\left(\sqrt{\ve^2+r^2}\right)~dr~\leq~{2C_3\ve^{\alpha}\over c_0}\cdot \big|{\bf v}_{11}(x)\big|\,.
\enda
\eeq
Putting together  (\ref{A1}), (\ref{Asm}),  (\ref{A11}), and (\ref{A1-2}), 
we can compare the sizes of the vectors $\bfv_{11}$ and $\bfv_{12}$ in 
(\ref{v11})--(\ref{v12}).   Indeed, the previous analysis shows that
\bel{A1-3} {|\bfv_{12}(x)|\over |\bfv_{11}(x)|}~\leq~
\ds{|A_1|+|A_2|+|A_3|\over \big|{\bf v}_{11}(x)\big|}~\leq~\ds{A_{11}+A_{12}+|A_2|+|A_3|\over \big|{\bf v}_{11}(x)\big|}~\to~0\qquad\mathrm{as}\qquad \ve~\to~0\,.
\eeq
\v
{\bf 3.} In a similar fashion we now compute
\bel{v2a}
{\bf v}_2(x)~=~{\bf v}_{21}(x)+{\bf v}_{22}(x),
\eeq
where 
\bel{v21}
\bega{rl}
{\bf v}_{21}(x)&=~\ds\int_{\ve^\alpha<|y|<\delta_0} \vp\Big(|x-y|
\Big)\cdot {x-y\over |x-y|}\, dy\\[4mm]
&=~\left(c_0\ds\int_{\ve^\alpha}^{\delta_0}r^{d-2} \vp\Big(\sqrt{\ve^2+r^2}
\Big)\cdot {\ve\over \sqrt{\ve^2+r^2}}\, dr\right)\, {\bf n}_x\,,
\enda
\eeq
and
 \bel{v22}\bega{rl}
 \bfv_{22}(x)&=~\ds\int_{\ve^\alpha<|y|<\delta_0}\Big[ \vp\bigl(|x-y-f(y) \bfn_x|\bigr)-
 \vp\bigl(|x-y|\bigr)\Big] \cdot
{x-y-f(y) \bfn_x\over|x-y-f(y) \bfn_x|}\, J(y)\, dy\\[4mm]
 &\qquad \ds+\int_{\ve^\alpha<|y|<\delta_0}\vp\bigl(|x-y|\bigr)\cdot
 \left\{{x-y-f(y) \bfn_x\over|x-y-f(y) \bfn_x|}-
{x-y\over|x-y|}\right\} \, J(y)\, dy
\\[4mm]
 &\qquad \ds+
 \int_{\ve^\alpha<|y|<\delta_0}\vp\bigl(|x-y|\bigr)\cdot
{x-y\over|x-y|}\, \bigl[J(y)-1\bigr]\, dy
\\[4mm]
 &=~B_1+B_2+B_3
 \,.\enda
\eeq
As in (\ref{A2})-(\ref{A3}), we have 
\bel{B23}
|B_2|+|B_3|~\leq~C_3\cdot \int_{\ve^{\alpha}}^{\delta_0}{r^{d}\over \sqrt{\ve^2+r^2}}\cdot \vp\left(\sqrt{\ve^2+r^2}\right)~dr~\leq~C_3\int_{\ve^{\alpha}}^{\delta_0}r^{d-1}\vp(\sqrt{\ve^2+r^2})~dr\,.
\eeq
Using again (\ref{ulb}) and the fact that $\vp$ is a decreasing function, we obtain
\bel{B1}
\bega{rl}
|B_1|&\leq~\ds C\cdot \Big[\int_{\ve^\alpha}^{\delta_0}
 r^{d-2} \, \left[\vp\left(\sqrt{\ve^2+(1-C_0\ve)r^2}\right)-\vp\Big(\sqrt{\ve^2+r^2}
\Big)\right] \, dr\, \\ [4mm]
&\ds\qquad\qquad\qquad+\int_{\ve^\alpha}^{\delta_0}
 r^{d-2} \, \left[\vp\Big(\sqrt{\ve^2+r^2}\Big)-\vp\Big(\sqrt{\ve^2+(1+C_0\ve)r^2+C_0^2r^4/4}
\Big)\right]dr\\[4mm]
&\doteq~B_{11}+B_{12}\,.
\enda
\eeq
As in (\ref{A1-1}), performing the variable change $s=\sqrt{1-C_0\ve}\cdot r$, we obtain
\bel{B1-1}
\bega{rl}\ds B_{11}&=~\ds C\cdot\int_{\ve^\alpha}^{\delta_0}
 r^{d-2} \, \left[\vp\left(\sqrt{\ve^2+(1-C_0\ve)r^2}\right)-\vp\Big(\sqrt{\ve^2+r^2}\Big)~dr\right]\\[4mm]
 &=~C\cdot\left[(1-C_0\ve)^{-{d-1\over 2}}\ds\int_{\ve^\alpha\sqrt{1-C_0\ve}}^{\delta_0\sqrt{1-C_0\ve}}s^{d-2}\vp\Big(\sqrt{\ve^2+s^2}\Big)~ds-\int_{\ve^{\alpha}}^{\delta_0}r^{d-2}\vp\Big(\sqrt{\ve^2+r^2}\Big)~dr\right]\\ [4mm]
 &\leq~\ds C_3\ve\int_{\ve^{\alpha}}^{\delta_0}r^{d-2}\vp\Big(\sqrt{\ve^2+r^2}\Big)~dr+C_3\, \int^{\ve^{\alpha}}_{\ve^{\alpha}\sqrt{1-C_0\ve}}r^{d-2}\vp\Big(\sqrt{\ve^2+r^2}\Big)~dr
\enda
\eeq
for a suitable constant $C_3$. Since $\vp$ is decreasing, for $\ve$ sufficiently small we have
\bel{est1}
\bega{rl}\ds
\int^{\ve^{\alpha}}_{0}r^{d-2}\vp\Big(\sqrt{\ve^2+r^2}\Big)~dr&\geq~\ds\int^{\ve^{\alpha}}_{\ve^{\alpha}(1-C_3\sqrt{\ve})}r^{d-2}\vp\Big(\sqrt{\ve^2+r^2}\Big)~dr\\ [4mm]
&\geq~\ds(1-C_3\sqrt{\ve})^{d-2}\cdot \int^{\ve^{\alpha}}_{\ve^{\alpha}(1-C_3\sqrt{\ve})}\ve^{d-2}\vp\Big(\sqrt{\ve^2+r^2}\Big)~dr\\ [4mm]
&\geq~\ds(1-C_3\sqrt{\ve})^{d-2}\cdot {C_3\ve^{\alpha+{1\over 2}}\over C_3\ve^{\alpha+1}}\cdot \int^{\ve^{\alpha}}_{\ve^{\alpha}(1-C_3\ve)}r^{d-2}\vp\Big(\sqrt{\ve^2+r^2}\Big)~dr\\ [4mm]
&\geq~\ds {1\over 2\sqrt \ve}\cdot \int^{\ve^{\alpha}}_{\ve^{\alpha}(1-C_3\ve)}r^{d-2}\vp\Big(\sqrt{\ve^2+r^2}\Big)~dr.
\enda
\eeq
In turn, this yields
\bel{B11}
\bega{rl}
 B_{11}&\leq~\ds C_3\ve\int_{\ve^{\alpha}}^{\delta_0}r^{d-2}\vp\Big(\sqrt{\ve^2+r^2}\Big)~dr+2C_3\sqrt{\ve}\int^{\ve^{\alpha}}_{0}r^{d-2}\vp\Big(\sqrt{\ve^2+r^2}\Big)~dr\\ [4mm]
 &\leq~\ds C_3\ve\int_{\ve^{\alpha}}^{\delta_0}r^{d-2}\vp\Big(\sqrt{\ve^2+r^2}\Big)~dr+2C_3\sqrt{\ve+\ve^{2\alpha-1}}\cdot\int_{0}^{\ve^{\alpha}}{\ve r^{d-2}\over \sqrt{\ve^2+r^2}}\cdot\vp\left(\sqrt{\ve^2+r^2}\right)~dr\\[4mm]
 &\leq~\ds C_3\ve\int_{\ve^{\alpha}}^{\delta_0}r^{d-2}\vp\Big(\sqrt{\ve^2+r^2}\Big)~dr+{4C_3\ve^{\alpha-{1\over 2}}\over c_0}\cdot  \big|{\bf v}_{11}(x)\big|\,.
 \enda
\eeq
Next, performing the change of variable $t=\sqrt{(1+C_0\ve)r^2+C_0^2r^4/4}$, we estimate 
\[
\bega{rl}
\ds\int_{\ve^\alpha}^{\delta_0}
 r^{d-2}\cdot\vp\Big(\sqrt{\ve^2+(1+C_0\ve)r^2+C_0^2r^4/4}\Big)\,dr&\geq~\ds\int^{\delta_0}_{\ve^{\alpha}(1+2C_0\sqrt{\ve})}{t^{d-2}\over  1+C_4(\ve+t^2)}\cdot \vp\left(\sqrt{\ve^2+t^2}\right)dt
 \enda
\]
for a suitable constant $C_4$. This implies that 
\[
\bega{rl}
\ds
B_{12}&\leq~C\cdot\left[\ds\int_{\ve^{\alpha}}^{\delta_0}r^{d-2}\vp\left(\sqrt{\ve^2+r^2}\right)dr-\int^{\delta_0}_{\ve^{\alpha}(1+2C_0\sqrt{\ve})}{r^{d-2}\over  1+C_4(\ve+r^2)}\cdot \vp\left(\sqrt{\ve^2+r^2}\right)dr\right]\\[4mm]
&\leq~C\cdot\left[\ds\int_{\ve^{\alpha}}^{\delta_0}\left(1-{1\over 1+C_4(\ve+r^2)}\right)r^{d-2}\vp\left(\sqrt{\ve^2+r^2}\right)dr+\int^{\ve^{\alpha}(1+2C_0\sqrt{\ve})}_{\ve^{\alpha}}r^{d-2}\vp\left(\sqrt{\ve^2+r^2}\right)dr\right]\\[4mm]
&\leq~\ds CC_4\int_{\ve^{\alpha}}^{\delta_0}(\ve+r^2)r^{d-2}\vp\left(\sqrt{\ve^2+r^2}\right)~dr+C\int^{\ve^{\alpha}(1+2C_0\sqrt{\ve})}_{\ve^{\alpha}}r^{d-2}\vp\left(\sqrt{\ve^2+r^2}\right)~dr\,.
\enda
\]
As in (\ref{est1}), one estimates
\[
\int^{\ve^{\alpha}(1+2C_0\sqrt{\ve})}_{\ve^{\alpha}}r^{d-2}\vp\left(\sqrt{\ve^2+r^2}\right)~\leq~2\sqrt{\ve}\int^{\ve^{\alpha}}_{0}r^{d-2}\vp\Big(\sqrt{\ve^2+r^2}\Big)dr
\]
for $\ve$ sufficiently small. Thus, as in (\ref{B1}), we obtain
\bel{B1-2}
\bega{rl}
B_{12}&\leq~\ds CC_4\int_{\ve^{\alpha}}^{\delta_0}(\ve+r^2)r^{d-2}\vp\left(\sqrt{\ve^2+r^2}\right)~dr+2C\sqrt{\ve}\int^{\ve^{\alpha}}_{0}r^{d-2}\vp\Big(\sqrt{\ve^2+r^2}\Big)dr\\[4mm]
&\leq~\ds CC_4\int_{\ve^{\alpha}}^{\delta_0}(\ve+r^2)r^{d-2}\vp\left(\sqrt{\ve^2+r^2}\right)~dr+4C\ve^{\alpha-{1\over 2}}\cdot \big|{\bf v}_{11}(x)\big|\,.
\enda
\eeq

Combining (\ref{v22}), (\ref{B23}), (\ref{B1}), (\ref{B11}), and (\ref{B1-2}), we finally obtain 
\bel{v22-1}
\left|{\bf v}_{22}(x)\right|~\leq~C_5\cdot \int_{\ve^{\alpha}}^{\delta_0}(\ve+r+r^2) r^{d-2}\vp\Big(\sqrt{\ve^2+r^2}\Big)~dr+C_5 \ve^{\alpha-{1\over 2}}\cdot \big|{\bf v}_{11}(x)\big|
\eeq
for some constant $C_5$.

\v{\bf 4.}   Finally, since in (\ref{v12}) the integral over $\Sigma_3$ involves functions 
which are uniformly bounded over  $\Sigma$, we have a trivial bound of the form
\bel{v3es}
|\bfv_3(x)|~\leq~C_6\,.\eeq

{\bf 5.} We now  compare the sizes of $\bfv_{22} (x)$ and $\bfv_3(x)$ with $\bfv_{11}(x)$. 
From (\ref{v11e}) it follows
\[
\bega{rl}
\bigl|\bfv_{11}(x)\bigr|~\geq~c_0\ds\int_{0}^{\ve^{\alpha}}r^{d}\cdot {\vp\left(\sqrt{\ve^2+r^2}\right)\over \sqrt{\ve^2+r^2}}~dr\,.
\enda
\]
Performing the change of variable $t=\sqrt{\ve^2+r^2}$ one obtains
\bel{v11-e1}
|{\bf v}_{11}(x)|~\geq ~c_0\ds\int_{\ve}^{\sqrt{\ve^2+\ve^{2\alpha}}}(t^2-\ve^2)^{d-1\over 2}\vp(t)~dt~\geq~c_0\ds\left({3\over 4}\right)^{d-1\over 2}\cdot \int_{2\ve}^{\ve^{\alpha}} r^{d-1}\vp(r)~dr.
\eeq
Recalling (\ref{v22-1}) and (\ref{v3es}), we obtain
\bel{ratio}
{ \bigl|{\bf v}_{22}(x)\bigr|+\bigl|{\bf v}_3(x)\bigr|
\over \big|{\bf v}_{11}(x)\big|}~\leq~C_7\cdot
{\ds  \int_{\ve^{\alpha}}^{\delta_0}r^{d-1}\vp\left(r\right)dr \over \ds \int_{2\ve}^{\ve^{\alpha}}r^{d-1}\vp(r)~dr}+C_5\cdot \ve^{\alpha-{1\over 2}}+{C_7\over |{\bf v}_{11}(x)|}\,,
\eeq
for some constant $C_7$. 

By (\ref{A1-3}) we already know that the ratio $|\bfv_{12}|/|\bfv_{11}|$
approaches zero as $\ve\to 0$.  Moreover, by (\ref{v11e}) and (\ref{v21}) one has
\[
{\bf v}_{21}(x)~=~C_{\ve}\, {\bf v}_{11}(x)\qquad\mathrm{for~some~constant~}C_{\ve}>0\,.
\]
Therefore, in view of (\ref{ratio}), we can conclude that 
(\ref{vperp}) holds true
provided that 
\bel{Ap1}
\lim_{\ve\to 0^+}~{\ds\int^{1}_{\ve^{\alpha}}r^{d-1}\vp(r)~dr\over \ds\int_{2\ve}^{\ve^{\alpha}}r^{d-1}\vp(r)~dr}~=~0\,.
\eeq
We show that (\ref{Ap1}) is  satisfied if $\vp$ satisfies the assumption {\bf (A2$'$)}. Indeed, let $\beta>{1\over2}$ be as in {\bf (A2$'$)} and choose ${1\over2}<\alpha<\beta$. Setting $\tilde{\ve}=2\ve$, we have $\ve^\alpha>\tilde{\ve}^\beta$ if $\ve$ is sufficiently small. Consequently,
\[
\ds\int_{\ve^{\alpha}}^1r^{d-1}\vp(r)dr~\leq~\ds\int_{\tilde{\ve}^{\beta}}^1r^{d-1}\vp(r)dr\qquad\mathrm{and}\qquad\ds\int_{2\ve}^{\ve^{\alpha}} r^{d-1}\vp(r)dr~\geq~\int_{\tilde{\ve}}^{\tilde{\ve}^{\beta}}r^{d-1}\vp(r)dr\,.
\]
On the other hand, as in the proof of Lemma 1, it holds
\[
\lim_{r\to 0+}~\int_{r}^{1}~s^{d-1}\vp(s)~ds~=~+\infty\,.
\]
Using L'Hopital's rule and the assumption ({\bf A2$'$}), we obtain 
\[
\lim_{\tilde{\ve}\to 0+}~{\ds\int_{\tilde{\ve}^{\beta}}^{1}r^{d-1}\vp(r)dr\over\ds \int_{\tilde{\ve}}^{1}r^{d-1}\vp(r)dr}~=~\lim_{\tilde{\ve}\to 0}~{\tilde{\ve}^{\beta(d-1)}\vp(\tilde{\ve}^{\beta})\cdot (\beta \tilde{\ve}^{\beta-1})\over \tilde{\ve}^{d-1}\cdot \vp(\tilde{\ve})}~=~\beta\cdot\lim_{\tilde{\ve}\to 0}~{\tilde{\ve}^{\beta d}\cdot \vp(\tilde{\ve}^{\beta})\over \tilde{\ve}^{d}\cdot \vp(\tilde{\ve})}~=~0\,.
\]
This implies 
\[
\lim_{\ve\to 0^+}~{\ds\int^{1}_{\ve^{\alpha}}r^{d-1}\vp(r)~dr\over \ds\int_{2\ve}^{\ve^{\alpha}}r^{d-1}\vp(r)~dr}~\leq~\lim_{\tilde{\ve}\to 0+}~{\ds\int_{\tilde{\ve}^{\beta}}^{1}r^{d-1}\vp(r)dr\over\ds \int_{\tilde{\ve}}^{\tilde{\ve}^{\beta}}r^{d-1}\vp(r)dr}~=~0,
\]
proving (\ref{Ap1}).
\endproof
\v
{\bf Corollary 2.} {\it Consider a family of compact $\C^2$ surfaces $\Sigma(t)$, 
continuously depending on $t\in [0,T]$, with uniformly bounded curvature.
Define the vector fields 
$\bfv(t,\cdot)$ as in (\ref{bvdef}).
Then for any $N,\ve>0$ there exists $\delta>0$ such that, for all $t\in [0,T]$,
\bel{vbp}\bega{rl}d(x, \Sigma(t))\,\leq\,\delta &\qquad\implies\qquad 
\bigl|\bfv(t,x)\bigr|\,\geq \,N,\\[2mm]
d(x, \Sigma(t))\,\leq\,\delta &\qquad\implies\qquad \ds
\left| {\bfv(t,x)\over \bigl|\bfv(t,x)\bigr|}- {x-\pi(t, x)\over \bigl|x-\pi(t, x)\bigr|}\right|\,\leq\,\ve.
\enda
\eeq
}

Indeed, the proof of Lemma~2 shows that the limits (\ref{vbig})-(\ref{vperp}) are uniformly valid over
a family of surfaces $\Sigma(t)$ with uniformly bounded curvature.

\v
\section{Proof of Theorem 3.} 
\label{sec:5}
\setcounter{equation}{0}  
Relying on Lemma 2, we can now give a proof of the convergence to the sweeping process, stated 
in (\ref{approx}).   We recall that this sweeping process keeps all trajectories inside a moving compact set $V(t)\subset \R^d$ with 
smooth boundary $\Sigma(t)$.   To fix the ideas, we assume that this set is defined in terms of a $\C^2$
function $\psi$, as in (\ref{Vdef})-(\ref{ppp}).
The argument relies on three main properties.
\begi
  \item[{\bf (P1)}] There exists a radius $\rho_0>0$ such that, if $t\in [0,T]$ and
 $d(x, \Sigma(t))\leq \rho_0$, then  the perpendicular 
 projection $\pi(t,x)$ of $x$ on $\Sigma(t)$ is well defined. In this case we denote by
 $\bfn(t,x)$  the unit normal vector  to $\Sigma(t)$  at   the  point $\pi(t,x)$, as in (\ref{ppro})
 \endi
 
\begi
 \item[{\bf (P2)}]  Setting
\bel{vtx}
\bfv(t,x)~\doteq~\int_{\Sigma(t)}\vp(|x-\xi|)\,  {x-\xi\over|x-\xi|}\, d\sigma(\xi),\eeq
for any $\delta>0$  the solution
 $t\mapsto x^{\delta}(t)$ to
 \bel{vdso}\dot x ~=~\delta \,\bfv(t,x),\qquad\quad x(0)= x_0\,.\eeq
 satisfies 
 $x^{\delta} (t)\in V(t)$, for every $x_0\in int(V(0))$ and $t\in [0,T]$.
 \endi 
 To see why this is true, assume $d(x^\delta(t), \Sigma(t))<\rho_0$
 and consider the unit normal vector
\bel{ppx}\bfn(t,x^\delta)~\doteq~
 {x^\delta(t)-\pi(t,x^\delta(t))\over |x^\delta(t)-\pi(t,x^\delta(t))|}.\eeq
Then
 $$
{d\over dt} \Big( d(x^\delta(t), \Sigma(t)) \Big)~\geq~\Big\langle
   \delta \bfv(t, x^\delta(t)),
\, \bfn(t, x^\delta)\Big\rangle - L_\Sigma\,,$$
where $L_\Sigma$ is a Lipschitz constant for the multifunction $t\mapsto \Sigma(t)$. 
By (\ref{vbig})-(\ref{vperp}), it follows that 
$$\la   \bfv(t, y),
\, \bfn(t,y)\ra ~\to~+\infty\qquad\hbox{as}\qquad d(y,\Sigma(t))\to 0.$$
Hence the distance $d(x^\delta(t), \Sigma(t))$ remains uniformly positive
in time, for every fixed $\delta>0$  and all initial points $x_0$ at a uniformly positive distance from $\Sigma(0)$.  
\v
Finally, by the properties (\ref{vplim}) of $\vp$ it follows
\begi
 \item[{\bf (P3)}]
 For every $0<\ve<{1\over 4}$, by choosing $0<\delta<\delta_0<\ve$ sufficiently small, for every $x\in V(t)$
 one has the implication
  \bel{alperp}0~<~d(x, \Sigma(t))~\leq~\delta_0
 \qquad\implies\qquad \left| {\bfv(t,x)\over |\bfv(t,x)|} - {x-\pi(t,x)\over |x-\pi(t,x)|} \right|~<~\ve.\eeq
 \bel{dvsmall}d(x, \Sigma(t))~>~\delta_0
 \qquad\implies\qquad \delta|\bfv(t,x)|~<~\ve.\eeq
\endi
Indeed, (\ref{alperp}) follows from Corollary~2 and the properties (\ref{ppp}) of the 
 function  $\psi$,  defining the boundary $\Sigma(t)$.  The implication (\ref{dvsmall}) trivially holds, 
 choosing  $\delta>0$ sufficiently small.
 \v
In the following, using the properties {\bf (P1)--(P3)}, we estimate the distance
between $x^\delta(t)$ and the solution $x(t,x_0)$ of the sweeping
process (\ref{sweep}). The proof will be given in 
several steps.
\v
{\bf 1.} For a given initial condition $x_0\in int V(0)$, let $t\mapsto x^\delta(t)$ be the 
solution to 
$$\dot x~=~\delta \bfv(t,x), \qquad \qquad x(0)= x_0\,,$$
and let $t\mapsto x(t)$ be the corresponding  solution to the sweeping process 
driven by the set $V(t)$.

For every $t\in[0,T]$ such that $d(x(t),\Sigma(t))\leq \rho_0/2$, let $\bfn(t)$ 
 be the unit normal vector  to $\Sigma(t)$  at   the  point $\pi(t,x(t))$.
By the regularity of $\Sigma(\cdot)$, we can extend $\bfn(\cdot)$ to a Lipschitz function defined 
on the entire time interval $[0,T]$. For simplicity, this extension will  
still be denoted by $t\mapsto\bfn(t)$.

We now split the difference as
\bel{wdef}w(t)~\doteq~x^\delta(t)-x(t)~=~w_1 (t)+ w_2(t), \eeq
where the vector $w_1(t)$ is parallel to $\bfn(t)$ while $w_2(t)$ 
is orthogonal to $\bfn(t)$.  Namely,
\bel{odec}
w_1(t)~\doteq~\theta(t) \bfn(t), \qquad \theta(t)~=~\la w(t), \bfn(t)\ra,
\qquad\quad w_2(t)~=~w(t)- w_1(t).\eeq
For future use, we observe that, if $d(x^\delta(t), \Sigma(t))\leq\rho_0$, then the unit normal vector
$\bfn(t,x^\delta)$ at (\ref{ppx}) 
is well defined and \bel{nnd}
|\bfn(t)-\bfn(t,x^\delta)|~\leq~C_\bfn |w(t)|,\eeq
for a suitable constant $C_\bfn$.

\begin{figure}[htbp]
\centering
\includegraphics[scale=0.5]{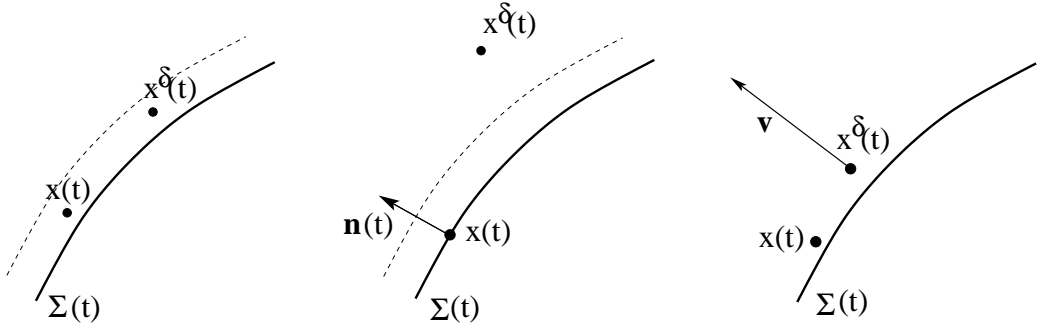}
    \caption{\small Left and center: the two different cases considered in the proof
    of Theorem~3, depending on the distance of $x^\delta(t)$
    from the boundary $\Sigma(t)$.  Right: if $d(x^\delta(t), \Sigma(t))<\delta_0$,
    then the speed $\delta\,\bfv(t,x^\delta)$ can be very large, and the same is true of 
    $\dot \theta$.   To handle this situation, we need to insert a weight function $W$
    in our estimates.}
\label{f:di14}
\end{figure}

\v
Our main goal is to show that $w(t)$ remains small.
This will be achieved by estimating the time derivatives $\dot w_1(t), \dot w_2(t)$, considering  two possible alternatives (see Fig.~\ref{f:di14}).

CASE 1: $d(x^\delta(t),\Sigma(t))\geq \delta_0\,$.

CASE 2: $d(x^\delta(t),\Sigma(t))<\delta_0\,$.   

We observe that, by the $\C^2$ regularity of the boundaries $\Sigma(t)= \partial V(t)$, there exists 
a constant $C_{\Sigma}$ such that
\bel{te0}
\la \bfn(t,x), \, y-x \ra~\geq~-C_{\Sigma}\cdot |y-x|^2\qquad\forall~~ t\in [0,T], ~~x\in\Sigma(t), ~~ y\in V(t)\,.
\eeq
In particular,
\bel{te} x(t)\in\Sigma(t)\qquad\implies\qquad 
\theta(t)~\geq~-C_{\Sigma}\, |w_2(t)|^2.\eeq

{\bf 2.} In this step we consider the Case 1: $d(x^\delta(t),\Sigma(t))\geq\delta_0\,$.

On the interval $[0,T]$,
 let $L_{{\bf n}}\geq 1$ be a Lipschitz constant for the map $t\mapsto \bfn(t)$, and let  $L_\Sigma$
 be a Lipschitz constant for the multifunction $t\mapsto \Sigma(t)$, w.r.t.~the Hausdorff distance.
 Observe that $L_\Sigma$ provides a common Lipschitz constant  for all trajectories $t\mapsto x(t)$ of 
 the sweeping process.

 We claim that
\bel{dtw}\ve + C_2| w_2(t)|~\geq~
\dot \theta(t)~\geq~\left\{\bega{cl}  -\ve-L_\bfn|w_2(t)|\ &\hbox{if}\quad x(t)\notin\Sigma(t)\\[1mm] -C_1 \ &\hbox{if}\quad x(t)\in\Sigma(t)\enda\right.,\qquad 
|\dot w_2(t)| ~\leq~\ve + C_2\,|w(t)|,\eeq
where $C_1=1+L_\Sigma+L_{{\bf n}}\cdot\max_{t\in [0,T]}\{\mathrm{diam}(V(t))\}$ and $C_2=2L_{\bf n}$. Indeed, recalling (\ref{dvsmall}) and (\ref{odec}), we have 
\bel{dt1}\bega{rl}
\dot{\theta}(t)&=~\langle \dot{w}(t),{\bf n}(t)\rangle+\langle w(t),\dot{{\bf n}}(t)\rangle\\[3mm]
&=~  \langle\dot{x}^{\delta}(t),{\bf n}(t)\rangle-\langle\dot{x}(t),{\bf n}(t)\rangle+\theta(t)\langle \bfn(t),\dot{{\bf n}}(t)\rangle+\langle w_2(t),\dot{{\bf n}}(t)\rangle\\[3mm]
&\geq~-|\dot{x}^{\delta}(t)|-|\dot{x}(t)|-|w_2(t)|\cdot |\dot{\bf n}(t)|\\[3mm]
&\geq~\left\{\bega{cl}  -\ve-L_\bfn|w_2(t)|\ &\hbox{if}\quad x(t)\notin\Sigma(t),\\[3mm] -\ve-L_\Sigma-L_\bfn|w_2(t)| \ &\hbox{if}\quad x(t)\in\Sigma(t).\enda\right.\enda 
\eeq
Moreover,
\bel{dt2}
\dot{\theta}(t)~\leq~\langle \dot{x}^{\delta}(t),{\bf n}(t)\rangle-\langle \dot{x}(t),{\bf n}(t)\rangle
+|w_2(t)|\cdot |\dot{\bf n}(t)|~\leq~\ve+L_{{\bf n}}\, |w_2(t)|.
\eeq
Together, (\ref{dt1})-(\ref{dt2}) yield
 the upper and lower bounds on $\dot\theta$  in (\ref{dtw}). 
 
Next, by (\ref{odec}) one has
\begin{eqnarray*}
|\dot{w}_2(t)|&=&\left|{d\over dt}~[w(t)-w_1(t)]\right|~=~\left|{d\over dt}\,\Big[w(t)-\langle w(t),{\bf n}(t)\rangle\, {\bf n}(t)\Big]\right|\\[3mm]
&=&\Big|\dot{w}(t)-\langle \dot{w}(t),{\bf n}(t)\rangle\, {\bf n}(t)-\langle w(t),\dot{{\bf n}}(t)\rangle\, {\bf n}(t)-\langle w(t),{\bf n}(t)\rangle\, \dot{{\bf n}}(t)\Big|\\[3mm]
&\leq&\Big|\dot{w}(t)-\langle \dot{w}(t),{\bf n}(t)\rangle\,{\bf n}(t)\Big|+2L_{{\bf n}}\, |w(t)|\,.
\end{eqnarray*}
Since $\dot x(t)$ is either zero or parallel to $\bfn(t)$, by
 (\ref{dvsmall}) it follows
\begin{eqnarray*}
\Big|\dot{w}(t)-\langle \dot{w}(t),{\bf n}(t)\rangle\, {\bf n}(t)\Big|&=&\Big|\delta\, {\bf v}(t,x^{\delta}(t))-\dot{x}(t)-\langle \delta\, {\bf v}(t,x^{\delta}(t))-\dot{x}(t),{\bf n}(t)\rangle\, {\bf n}(t)\Big|\\[3mm]
&=&|\delta\, {\bf v}(t,x^{\delta}(t))-\langle \delta\, {\bf v}(t,x^{\delta}(t)),{\bf n}(t)\rangle\, {\bf n}(t)|\\[3mm]
&\leq&|\delta\, {\bf v}(t,x^{\delta}(t))|~\leq~\ve\,.
\end{eqnarray*}
This implies the second inequality in (\ref{dtw}). In particular, if $x(t)\in int(V(t))$  then $\dot x(t)=0$ and
\bel{ot1}
{d\over dt}|w(t)|~\leq~3(\ve+{\bf L}_\bfn\, |w(t)|)\,.
\eeq
\v
{\bf 3.} In this step we consider  Case 2: $d(x^\delta(t),\Sigma(t))<\delta_0\,$.
As long as
\bel{wsma} |w(t)|~\leq~{1\over 4C_{\bf n}},\eeq
we claim that 
\bel{tb1}
\theta(t)~\leq~\delta_0 + C_3\bigl( \delta_0 + |w(t)|\bigr)^2,\eeq
\bel{td1}
\dot\theta(t)~\geq~\left\{\bega{cl}  -L_\bfn|w_2(t)|\ &\hbox{if}\quad x(t)\notin\Sigma(t),\\[3mm] -C_1 \ &\hbox{if}\quad x(t)\in\Sigma(t),\enda\right.\qquad |\dot w_2(t)| ~\leq~
C_4\bigl( \ve+|w(t)|\bigr) \bigl(1+ |\dot \theta(t)|\bigr),\eeq
for some constants $C_3, C_4$.

Notice that, if $d(x(t), \Sigma(t))\geq\rho_0/2$, we then have $|w(t)|\geq \rho_0/4$, because without loss of generality we can assume $\delta_0<\ve<\rho_0/4$. 
In this case the estimate (\ref{tb1}) is trivially satisfied, by choosing a constant 
$C_3$ large enough.

In the following, we thus assume $d(x(t), \Sigma(t))\geq\rho_0/2$, so that the projection $\pi(t,x(t))$ is well
defined. 
Recalling (\ref{wdef})--(\ref{te0}) we obtain
\begin{eqnarray*}
\theta(t)&=& \left\langle x^{\delta}(t)-x(t),{\bf n}(t)\right\rangle~=~ \left\langle x^{\delta}(t)-\pi(t,x^{\delta}),{\bf n}(t)\right\rangle+ \left\langle \pi(t,x^{\delta})-x(t),{\bf n}(t)\right\rangle\\[3mm]
&\leq& \left|x^{\delta}(t)-\pi(t,x^{\delta})\right|+\left\langle \pi(t,x^{\delta})-x(t),{\bf n}(t)\right\rangle
\\[3mm]
&\leq&\delta_0+\left\langle \pi(t,x^{\delta})-x(t)
,\, {\bf n}(t)-{\bf n}(t,x^\delta)\right\rangle+\left\langle\pi(t,x^{\delta})-x(t),\, 
{\bf n}(t,x^{\delta})\right\rangle
\\[3mm]
&\leq&\delta_0+C_{{\bf n}}\, |w(t)|\, \left|\pi(t,x^{\delta})-x(t)\right|+
C_\Sigma \,\left|\pi(t,x^{\delta})-x(t)\right|^2
\\[3mm]
&\leq&\delta_0+C_{{\bf n}}\, |w(t)|\, \left(\delta_0+|w(t)|\right)+
C_{\Sigma}\, \left(\delta_0+|w(t)|\right)^2\\[3mm]
&\leq&\delta_0+C_3\left(\delta_0+|w(t)|\right)^2
\end{eqnarray*}
for a suitable constant $C_3$.
This implies (\ref{tb1}). 

Using (\ref{alperp}), the time derivative $\dot \theta$ can be estimated as
\bel{dt6}\bega{rl}
\dot{\theta}(t)&=~\langle \dot{w}(t),{\bf n}(t)\rangle+\langle w(t),\dot{{\bf n}}(t)\rangle~\geq~ \langle \dot{x}^\delta(t),{\bf n}(t)\rangle- |\dot{x}(t)|-|w_2(t)|\, |\dot{\bf n}(t)|
\\[3mm]
&=~\delta |\bfv(t,x^\delta)|\left\langle \Big({\bfv(t,x^\delta)\over |\bfv(t,x^\delta)|}-\bfn(t,x^\delta)\Big)
+\Big(\bfn(t,x^\delta)-\bfn(t)\Big)+\bfn(t),\,{\bf n}(t)\right\rangle\\[3mm]
&\qquad- |\dot{x}(t)|-|w_2(t)|\, |\dot{\bf n}(t)|
\\[3mm]
&\geq~\left\{\bega{ll}\delta |\bfv(t,x^\delta)|\Big(-\ve-C_\bfn |w(t)|+1\Big)-L_\bfn\, |w_2(t)|\ &\hbox{if}\quad x(t)\notin\Sigma(t),\\[3mm]\delta |\bfv(t,x^\delta)|\Big(-\ve-C_\bfn |w(t)|+1\Big)-L_\Sigma-L_\bfn\, |w_2(t)|\ &\hbox{if}\quad x(t)\in\Sigma(t).\enda\right.
\enda\eeq
As long as (\ref{wsma}) holds,  we have $1-\ve-C_\bfn |w(t)|\geq {1\over 2}$. This already yields
the first inequality in (\ref{td1}).
{}From (\ref{dt6}) we also deduce
\bel{dbf}
 |\delta\bfv(t,x^\delta)|~\leq~2\bigl(C_1+ |\dot \theta(t)|\bigr).
 \eeq
In turn, this yields
\begin{eqnarray*}
|\dot{w}_2(t)|&\leq&\left|\delta {\bf v}(t,x^{\delta})-\langle \delta\, {\bf v}(t,x^{\delta}),{\bf n}(t)\rangle\,
{\bf n}(t)\right|+2L_{{\bf n}}\, |w(t)|
\\[3mm]
&\leq&|\delta\bfv(t,x^\delta)|\, \left| {\bfv(t,x^\delta)\over |\bfv(t,x^\delta)|}-\bfn(t)\right|
+2L_{{\bf n}}\, |w(t)|\\[3mm]
&\leq&|\delta\bfv(t,x^\delta)|\,\left\{\left| {\bfv(t,x^\delta)\over |\bfv(t,x^\delta)|}-\bfn(t,x^{\delta})\right|+\left|{\bf n}(t,x^{\delta})-{\bf n}(t)\right|\right\}+2L_{{\bf n}}\, |w(t)|\\[3mm]
&\leq&|\delta\bfv(t,x^\delta)|\, (\ve+C_{\bf n} |w(t)|)+2L_{\bf n}\, |w(t)|\\[3mm]
&\leq&2\bigl(C_1+ |\dot \theta(t)|\bigr)\, (\ve+C_{\bf n} |w(t)|)+2L_{\bf n}\, |w(t)|,
\end{eqnarray*}
establishing the second inequality in (\ref{td1})\,.
\v
{\bf 4.} 
In this step we prove that there exists a constant $C_5>0$ such that, for every $t\in[0,T]$,
\bel{tb2}|\theta(t)|~\leq~ C_5(\ve+\Tilde{w}_2(t))\,,
\eeq
where 
\bel{tdww}
\Tilde{w}_2(t)~\doteq~\max_{s\in [0,t]}~|w_2(s)|\qquad\forall t\in [0,T].
\eeq
Notice that, if $x(t)\in \Sigma(t)$, then (\ref{te}) applies. Let us assume that $x(t)\notin \Sigma(t)$ and
define the time
$$t_0~\doteq~\inf\bigl\{s\in [0,t]~;~x(r)\notin \Sigma(r)\quad\forall r\in [s,t]\bigr\}\,.$$
From (\ref{dtw}) and (\ref{td1}) it follows 
\[
\dot{\theta}(s)~\geq~-\ve-L_{\bf n} |w_2(s)|)~\geq~-\ve-L_{\bf n}\Tilde{w}_2(s)\qquad\hbox{for a.e.}~ s\in [t_0,t]\,.
\]
This yields
\begin{eqnarray*}
\theta(t)&=&\theta(t_0)+\int_{t_0}^t\dot{\theta}(s)\,ds~\geq~\theta(t_0)-\int_{t_0}^{t}\ve+L_{\bf n}\Tilde{w}_2(s)~ds\\[3mm]
&\geq&-C_{\Sigma}\, |w_2(t_0)|^2-T(\ve+L_\bfn\Tilde{w}_2(t))~\geq~-C_{\Sigma}\, \Tilde{w}_2^2(t)-T(\ve+L_\bfn\Tilde{w}_2(t))\,.
\end{eqnarray*}
As long as $\Tilde w_2(t)<1$, recalling that $L_\bfn\geq 1$ we have 
\[
\theta(t)~\geq~-\left(C_{\Sigma}+TL_\bfn\right)\cdot(\ve+\Tilde{w}_2(t))\,.
\]
Thus, to obtain (\ref{tb2}), one only needs to  consider the case where $\theta(t)>0$. Observe that, if $d(x^\delta(t),\Sigma(t))\leq\delta_0$, as long as
\bel{wsma1}
\ve+|w(t)|~\leq~{1\over 2C_3}\,,
\eeq
from (\ref{tb1}) it follows
\bel{thb5}\theta(t)~\leq ~\Tilde{C}_5(\ve+\Tilde{w}_2(t))\eeq
for some constant $\Tilde{C}_5$.  

We claim that  (\ref{tb2}) holds, for a suitable constant $C_5$.
Indeed, consider the time
$$t_1~\doteq~\inf\Big\{s\in [0,t)\,; ~\theta(r)>0\quad\hbox{and}\quad d(x^\delta(r),\Sigma(r))> \delta_0\quad\forall r\in (s,t)\Big\}\,.$$
We then have $\theta(t_1)\leq \Tilde{C}_5(\ve+\Tilde{w}_2(t))$.  Therefore, by (\ref{thb5}) and (\ref{dtw}) 
it follows
$$\bega{rl}
\theta(t)&=~\ds\theta(t_1)+\int_{t_1}^t\dot{\theta}(s)\,ds~\leq~\Tilde{C}_5(\ve+\Tilde{w}_2(t_1))+\int_{t_1}^{t}
\Big(\ve + C_2| w_2(s)|\Big)\,ds\\[3mm]
&\leq~\ds\Tilde{C}_5(\ve+\Tilde{w}_2(t))+\int_{t_1}^{t}\Big(\ve + C_2\Tilde  w_2(t)
\Big)\,ds~\leq~C_5 (\ve+\Tilde{w}_2(t)).
\enda $$
\v
{\bf 5.} In the following we shall assume $|w(t)|\leq \rho_0/3$ for all 
$t\in [0,T]$.   In this case, when $d(x(t), \Sigma(t))>\rho_0$ 
we have $\dot x(t)=0$, $|\dot x^\delta(t)|<\ve$ and the estimates are trivial.
Without loss of generality, we can thus assume that the 
normal vectors $\bfn(t)$ and $\bfn(t, x^\delta)$ are well defined.

For a suitable constant $\kappa$ (to be determined later), define
the weight
\bel{W} W(t)~\doteq~\left\{\bega{cl} 
\exp\Big\{ -\kappa \, d(x^\delta(t), \Sigma(t))\Big\}\qquad &\hbox{if}\quad d(x^\delta(t),\Sigma(t))\leq \delta_0\,,\\[4mm]
\exp\{ -\kappa \delta_0\}\qquad &\hbox{if}\quad d(x^\delta(t),\Sigma(t))\geq \delta_0\,.
\enda\right.\eeq
We now analyze how the weighted distance
$$\Lambda(t)~\doteq~|\theta(t)|+ W(t)\, \Tilde{w}_2(t)$$
changes in time.   The heart of the matter is to provide a bound on $\Tilde{w}_2$.
Indeed, by (\ref{tb2}) the component $w_1(t)=\theta(t) \bfn(t)$ 
can be bounded
in terms of $\Tilde{w}_2(t)$.
\v
{\bf 6.} 
At any point $t\in [0,T]$ where  ${w}_2(\cdot)$ is differentiable, by the definition  of $\Tilde w_2$ it follows
\bel{bdtw}
0~\leq~{d\over dt}~\Tilde{w}_2(t)~\leq~|\dot{w}_2(t)|.
\eeq
%
We first consider Case 1, where $d(x^\delta(t), \Sigma(t))\geq\delta_0$. 
By (\ref{dtw}) and (\ref{tb2}), we have that $\dot{W}(t)=0$ and
\[|\dot w_2(t)| ~\leq~\ve + C_2\,|w(t)|~\leq~C_2\,(\ve + |\theta(t)|+\Tilde{w}_2(t))~\leq~C_2\,(1+C_5)(\ve+\Tilde{w}_2(t))\,.\]
Therefore, (\ref{bdtw}) yields
\begin{eqnarray*}
{d\over dt}~\bigl(W(t)\Tilde{w}_2(t)\bigr)&=&W(t)\,{d\over dt}\Tilde{w}_2(t)~\leq~W(t)\cdot |\dot{w}_2(t)|~\leq~C_2\,(1+C_5)W(t)(\ve+\Tilde{w}_2(t))\,.
\end{eqnarray*}

On the other hand, 
in the case where
$d(x^\delta(t), \Sigma(t))<\delta_0$, by (\ref{td1}) and (\ref{tb2}) one has
\[
\bega{rl}
\ds \left|\dot{w}_2(t)\right|&\leq~C_4(\ve+|w(t)|) (1+|\dot\theta(t)|)\\[3mm]
&\leq~C_4(\ve+|\theta(t)|+\Tilde{w}_2(t)) (1+|\dot\theta(t)|)~\leq~(C_4+C_4C_5)(\ve+\Tilde{w}_2(t)) (1+|\dot\theta(t)|)\\[3mm]
&\leq~\left\{ \bega{cl} (C_4+C_4C_5)(1+C_1)\,(\ve+\Tilde{w}_2(t)) \quad &\hbox{if}\quad\dot\theta(t)\leq 0,
\\[3mm]
(C_4+C_4C_5) (\ve+\Tilde{w}_2(t)) (1+\dot\theta(t))\quad &\hbox{if}\quad\dot\theta(t)> 0\,,
\enda\right.\\
\quad\\
&\leq~\left\{ \bega{cl} C_6(\ve+\Tilde{w}_2(t)) \quad &\hbox{if}\quad\dot\theta(t)\leq 0,
\\[4mm]
C_6 (\ve+\Tilde{w}_2(t)) (1+\dot\theta(t))\quad &\hbox{if}\quad\dot\theta(t)> 0\,,
\enda\right.
\enda
\]
%
%
%
By (\ref{bdtw}), for a.e.~$t\in [0,T]$ one has
\bel{dtw2}
{d\over dt}\Tilde{w}_2(t)~\leq~\left\{ \bega{cl} C_6(\ve+\Tilde{w}_2(t)) \quad &\hbox{if}\quad\dot\theta(t)\leq 0,
\\[4mm]
C_6 (\ve+\Tilde{w}_2(t)) (1+\dot\theta(t))\quad &\hbox{if}\quad\dot\theta(t)> 0\,.
\enda\right.
\eeq

As long as $\ve+|w(t)|\leq{1\over 2C_1(2C_{\bf n}+C_4)}$,  from (\ref{nnd}), (\ref{dt1}), and (\ref{dbf}), 
it follows
\[
\bega{l}\ds
{d\over dt}~d(x^\delta(t), \Sigma(t))~\geq~\Big\langle\delta \bfv(t, x^\delta(t)),\, \bfn(t, x^\delta)\Big\rangle-L_{\Sigma}\\[4mm]
\qquad \geq~\Big\langle\delta \bfv(t, x^\delta(t)),\, \bfn(t)\Big\rangle- C_{\bf n}|w(t)|\cdot |\delta \bfv(t, x^\delta(t))|-L_{\Sigma}\\[4mm]
\qquad \geq~\Big\langle\dot{w}_1(t),\, \bfn(t)\Big\rangle+\Big\langle\dot{w}_2(t),\, \bfn(t)\Big\rangle+\Big\langle\dot{x}(t),\, \bfn(t)\Big\rangle- C_{\bf n}|w(t)|\cdot |\delta \bfv(t, x^\delta(t))|-L_{\Sigma}\\[4mm]
\qquad \geq~\dot{\theta}(t)-|\dot{w}_2(t)|-C_{\bf n}|w(t)|\cdot |\delta \bfv(t, x^\delta(t))|-L_{\Sigma}\\[4mm]
\qquad \geq~\dot{\theta}(t)-C_4(\ve+|w(t)|)\cdot \bigl(1+ |\dot \theta(t)|\bigr)-2C_{\bf n}|w(t)|\cdot \bigl(C_1+ |\dot \theta(t)|\bigr)-L_{\Sigma}\\[4mm]
\qquad \geq~\dot{\theta}(t)-(2C_{\bf n}+C_4)(\ve+|w(t)|)\cdot \bigl(C_1+ |\dot \theta(t)|\bigr)-L_{\Sigma}~\geq~{1\over 2}\dot{\theta}(t)-C_7
\enda
\]
for some constant $C_7>0$. Inserting the weight, we now estimate 
\begin{eqnarray*}
{d\over dt}\Big(W(t)\Tilde{w}_2(t)\Big)&=&-\kappa W(t)\Tilde{w}_2(t)\cdot {d\over dt}d(x^\delta(t), \Sigma(t))+W(t)\cdot {d\over dt} {\Tilde{w}}_2(t)\\[4mm]
&\leq&-{\kappa\over 2}W(t)\Tilde{w}_2(t)\dot{\theta}(t)+\kappa C_7 W(t)\Tilde{w}_2(t)+W(t)\cdot {d\over dt}{\Tilde{w}}_2(t)\,.
\end{eqnarray*}
Two cases can occur:
\begin{itemize}
\item If $\dot\theta(t)\leq 0$, then (\ref{dtw2}) and (\ref{td1}) yield
\[
{d\over dt}\Big(W(t)\Tilde{w}_2(t)\Big)~\leq~C_8 \, W(t)\, (\ve + \Tilde{w}_2(t))
\]
for some constant $C_8$\,.
\item If $\dot\theta(t)> 0$, then (\ref{dtw2}) yields
$$\bega{l}\ds
{d\over dt}~\left(W(t)\Tilde{w}_2(t)\right)~\leq~-{\kappa\over 2}W(t)\Tilde{w}_2(t)\dot{\theta}(t)+W(t)\cdot{d\over dt}{\Tilde{w}}_2(t)+\kappa C_7 W(t)\Tilde{w}_2(t)\\[4mm]
\qquad\ds \leq~-{\kappa\over 2}W(t)\Tilde{w}_2(t)\dot{\theta}(t)+C_6 W(t)(\ve+\Tilde{w}_2(t)) (1+\dot\theta(t))+\kappa C_7 W(t)\Tilde{w}_2(t)\\[4mm]
\qquad\ds \leq~W(t)\cdot\left[\left(-{\kappa\over 2}\Tilde{w}_2(t)+
C_6\cdot (\ve+\Tilde{w}_2(t))\right)\cdot \dot{\theta}(t)+(C_6+\kappa C_7)(\ve+\Tilde{w}_2(t))\right]\,.
\enda
$$
We now choose the constant $\kappa$ in (\ref{W}) so that 
$${\kappa\over 2}~\geq~2 C_6\,.$$
In this case, either $\Tilde w_2(t)<\ve$, or else
\bel{dtw5}\bega{rl}
\ds{d\over dt}\Big( W(t)\, \Tilde{w}_2(t)\Big)&\ds\leq~W(t)\cdot \left[\left( -{\kappa\over 2} + 2C_6\right) 
\dot\theta(t)\Tilde{w}_2(t) + (C_6+\kappa C_7)(\ve+\Tilde{w}_2(t))\right]\\[4mm]
&\ds\leq~ W(t)\,(C_6+\kappa C_7)(\ve+\Tilde{w}_2(t))\,.\enda\eeq
\end{itemize}
Combining both Cases 1 and 2, we obtain that either $\Tilde{w}_2(t)\leq \ve$ or else 
\bel{dW}{d\over dt}\Big( W(t)\, \Tilde{w}_2(t)\Big)~\leq~C_9 \,W(t)\,\Tilde{w}_2(t),\eeq
provided that 
\bel{wsma2}
\ve+|\theta(t)|+\Tilde{w}_2(t)~\leq~\min\left\{{1\over 4},{1\over 4C_{\bf n}},{1\over 2C_3},{1\over 2C_1(2C_{\bf n}+C_4)},{\rho_0\over 3}\right\}\,.
\eeq
%
%
\v
{\bf 7.} To complete the argument, consider the time
\[
\bar{t}~\doteq~\sup\,\Big\{\tau\in [0,T]\,;~(\ref{wsma2})~\mathrm{holds~for~all}~t\in [0,\tau]\Big\}\,.
\] 
Since $\Tilde{w}_2(t)$ is continuous and non-decreasing, there exists $t_{\ve}\in [0,\bar{t}]$ such that 
\[
\left\{ \bega{cl} \Tilde{w}_2(t)~\leq~\ve\qquad\forall t\in [0,t_{\ve}],
\\[4mm]
\Tilde{w}_2(t)~>~\ve\qquad\forall t\in \,]t_{\ve},\bar{t}]\,.
\enda\right.
\]
Hence (\ref{dW}) implies  
\[
W(t)\Tilde{w}_2(t)~\leq~e^{C_9 (t-t_{\ve})}W(t_{\ve})\Tilde{w}_2(t_{\ve})\qquad\forall t_{\ve}\leq t\leq \bar{t}\,.
\]
Since $e^{-\kappa\delta_0}\leq W(t)\leq1$, we have
\[
\Tilde{w}_2(t)~\leq~\exp\left({C_9T+\kappa}\cdot \delta_0\right)\cdot \ve\qquad\forall t\in [0,\bar{t}]\,.
\]
Recalling (\ref{tb2}), we obtain
\[
|\theta(t)|~\leq~C_5\,\left[\exp\left({C_9T+\kappa}\cdot \delta_0\right)+1\right]\cdot \ve\qquad\forall t \in [0,\bar{t}]\,.
\]
This yields 
\[
\ve+|\theta(t)|+\Tilde{w}_2(t)~\leq~C_{10}\,\ve\qquad\forall t\in [0,\bar{t}],
\]
where 
$$C_{10}~=~(1+C_5)\left[\exp\left({C_9T+\kappa}\cdot \delta_0\right)+1\right]\,.$$
Therefore, for any $\ve>0$ such that $$\ve~<~{1\over C_{10}}\min\left\{{1\over 4},\,{1\over 4C_{\bf n}},\,
{1\over 2C_3},\,{1\over 2C_1(2C_{\bf n}+C_4)},\,{\rho_0\over 3}\right\},$$ we conclude that 
\bel{swap2}\bar{t}=T,\qquad\qquad 
|w(t)|~\leq~|\theta(t)|+\Tilde{w}_2(t)~\leq~C_{10}\, \ve\qquad\forall t\in [0,T]\,.
\eeq
\v
{\bf 8.} The previous analysis has shown that, by choosing $\delta>0$ small enough, the sweeping process 
can be arbitrarily well approximated by 
the evolution generated by the vector field $\delta\bfv(t,x)$.
Repeating the argument in step {\bf 3} of the proof of Theorem 2,
we now construct  
a control function  $t\mapsto \xi(t)$ such that trajectories of the ODE
$$\dot x(t)~=~\vp(|x-\xi(t)|){x-\xi(t)\over |x-\xi(t)|}$$
 approximate the trajectories  of $\dot x = \delta\,\bfv(t,x)$, uniformly for $t\in [0,T]$ and
for all initial data in the compact set $\Omega_0\subset \,\R^d\setminus\Sigma(0)$.
This completes the proof.
\endproof
\v
{\bf Acknowledgments.} The research by K.~T.~Nguyen was partially supported by a grant from the Simons Foundation/SFARI (521811, NTK).

\end{document}